\documentclass[12pt]{article}
\usepackage{amsmath,pslatex,amssymb,amsfonts,latexsym,amstext,amsthm,enumerate,color}
\usepackage{setspace,booktabs,multirow,colortbl,float}
\usepackage{todonotes}

\newtheorem{theorem}{Theorem}
\newtheorem{lemma}[theorem]{Lemma}
\newtheorem{corollary}[theorem]{Corollary}
\newtheorem{conjecture}[theorem]{Conjecture}

\newtheorem{definition}[theorem]{Definition}

\newtheorem{proposition}[theorem]{Proposition}

\newcommand\Zp{\makebox[.75em]{$1$}}
\newcommand\Zm{\makebox[.75em]{$-$}}
\newcommand\Zz{\makebox[.75em]{$0$}}

\newcommand\Za{\makebox[.75em]{$a$}}
\newcommand\ZA{\makebox[.75em]{$\overline{a}$}}
\newcommand\Zb{\makebox[.75em]{$b$}}
\newcommand\ZB{\makebox[.75em]{$\overline{b}$}}

{\begin{list}{}%
         {\setlength{\leftmargin}{#1}}%
         \item[]%
}
{\end{list}}

\definecolor{Gray}{gray}{0.9}
\newcommand{\ra}[1]{\renewcommand{\arraystretch}{#1}}
\ra{1.2}

\newcounter{mycomment}

\begin{document}

\title{Mutually Unbiased Weighing Matrices}

\author{
D. Best\footnote{Supported by NSERC CGS-M and Alberta Innovates -- Technology Futures.}, ~~ H. Kharaghani\footnote{Supported by an NSERC Discovery Grant. Corresponding author.}\\
Department of Mathematics and Computer Science\\University of Lethbridge, AB T1K 3M4, Canada \\ darcy.best@uleth.ca, kharaghani@uleth.ca
\and
H. Ramp\\
Department of Physics\\University of Alberta, AB T6G 2R3, Canada \\ ramp@ualberta.ca
}

\maketitle

\begin{abstract}
Inspired by the many applications of mutually unbiased Hadamard matrices, we study mutually unbiased weighing matrices. These matrices are studied for small orders and weights in both the real and complex setting. Our results make use of and examine the sharpness of a very important existing upper bound for the number of mutually unbiased weighing matrices.
\end{abstract}



\section{Introduction} \label{sec:intro}

A {\it unit weighing matrix}, $W=[w_{ij}]$, with order $n$ and weight $p$, denoted $UW(n,p)$, is an $n \times n$ matrix with $|w_{ij}|$ in $\left\{0,1\right\}$ and  $WW^*=pI_n$, where $W^*=[{\overline{w}}_{ji}]$ is the the usual conjugate transpose of $W$. This implies that the rows of $W$ are mutually orthogonal under the standard inner product in $\mathbb{C}^n$ and contain exactly $p$ nonzero entries in each row and column. When $n = w$ (i.e., no zeroes in the matrix), $W$ is a {\it Hadamard matrix}. A {\it real weighing matrix} is the one with $w_{ij}$ in $\left\{0,\pm1\right\}$. Real weighing matrices have been well studied for small weights (see \cite{inequiv-weigh}) and large weights (see \cite{Craigen_1995}). This article contains results for weighing matrices in both the real and complex setting. Motivated by the applications of real weighing matrices, we have studied unit weighing matrices in \cite{unit-weigh13}. Our aim in this paper is to complement the work in \cite{unit-weigh13}.

Two unit weighing matrices $UW(n,p)$, $H$ and $K$, are {\it unbiased} if $HK^*=\sqrt{p}L$, where $L$ is a unit weighing matrix $UW(n,p)$. A set of pairwise unbiased unit weighing matrices are called {\it mutually unbiased unit weighing matrices}. In the special case of $n = w$, these are termed mutually unbiased Hadamard matrices (MUHM), which are of great interest to people working in areas related to the quantum information theory and as such, there is extensive literature on these matrices. We refer the reader to the most comprehensive survey paper \cite {durt-muhm} on MUHM. Mutually unbiased unit weighing matrices have also seen some application in quantum information science, specifically in the context of zero-error classical communication. \cite {Leung}

In \cite{unit-weigh13}, we concerned ourselves with the existence of certain unit weighing matrices; here, we are concerned about how many pairwise unbiased unit weighing matrices there are. In the general unimodular case, we lose a lot of structure that can be found in the real case (see Lemma \ref{lem:perf-sqr} for one such example), which makes it very challenging to locate complete sets.

If the entries of matrices in a set of mutually unbiased unit weighing matrices are limited to certain roots of unity, then a bound similar to Lemma \ref{lem:perf-sqr} is found (ex., see \cite{much10}), but very few concrete bounds exist in general. Section \ref{sec:restrictions} will deal with the unit weighing matrices in general by giving the few known upper bounds and lower bounds on the size of these sets.

In Section \ref{sec:main}, we will outline some of our computer searches for small orders of real weighing matrices. As an extension to mutually unbiased unit weighing matrices, we will examine sets of Hadamard matrices whose pairwise products satisfy specific conditions in section \ref{sec:weakly-unbiased}.

\section{General Restrictions}\label{sec:restrictions}

We will start off with a very well-known result (see \cite{much10}).

\begin{lemma}\label{lem:perf-sqr}
 Let $H$ and $K$ be real unbiased weighing matrices with order $n$ and weight $w$, then $w$ must be a perfect square.
 \begin{proof}
  Since both $H$ and $K$ are integer matrices, $HK^T=\sqrt{w}L$ must be an integer matrix as well.
 \end{proof}
\end{lemma}

\begin{lemma}\label{lem:single-direct-sum}
 Let $\mathcal{W} = \left\{W_1,\cdots,W_k\right\}$ be a set of mutually unbiased weighing matrices of order $m$ with weight $w$ and $\mathcal{X} = \left\{X_1,\cdots,X_l\right\}$ be a set of mutually unbiased weighing matrices of order $n$ with weight $w$. Then there exist $p=\min(k,l)$ mutually unbiased weighing matrices of order $m+n$ and weight $w$.
 \begin{proof}
  The set $\left\{W_1 \oplus X_1, W_2 \oplus X_2, \cdots, W_p \oplus X_p\right\}$ gives the desired result, where $\oplus$ denotes the standard direct sum of matrices (i.e., $A \oplus B = \text{diag}(A,B)$).
 \end{proof}
\end{lemma}

\begin{theorem}\label{th:direct-sum}
 Let $\left\{\mathcal{W}_1,\cdots,\mathcal{W}_k\right\}$ be a collection of sets of mutually unbiased weighing matrices of order $n_i$ and weight $w$. Then there are $$\min_{1 \leq i \leq k}\left(\left|\mathcal{W}_i\right|\right)$$ mutually unbiased weighing matrices of order $\sum_{i=1}^{k} n_i$ and weight $w$.
 \begin{proof}
 The case where $k=1$ is trivially true. Now assume the property holds for a collection of size $k-1 \geq 1$. Consider a collection with $k$ elements. By applying Lemma \ref{lem:single-direct-sum} to $\mathcal{W}_1$ and $\mathcal{W}_2$, we know there exists a collection of mutually unbiased weighing matrices of order $n_1+n_2$ and weight $w$ with $min\left(\left|\mathcal{W}_1\right|,\left|\mathcal{W}_2\right|\right)$ elements (we shall call this collection $\mathcal{X}$). By applying the induction hypothesis to 
$\left\{\mathcal{X},\mathcal{W}_3,\cdots,\mathcal{W}_{k}\right\}$, 
we have that there are 
$$\min\left(\left|\mathcal{X}\right|,\left|\mathcal{W}_3\right|, \cdots, \left|\mathcal{W}_{k}\right|\right)
=\min\left(\min\left(\left|\mathcal{W}_1\right|,\left|\mathcal{W}_2\right|\right),\left|\mathcal{W}_{3}\right|, \cdots, \left|\mathcal{W}_{k}\right|\right)
= \min_{1 \leq i \leq k}\left(\left|\mathcal{W}_i\right|\right)$$ mutually unbiased weighing matrices of order $(n_1+n_{2})+\sum_{i=3}^{k} n_i=\sum_{i=1}^{k} n_i$ and weight $w$.

 \end{proof}
\end{theorem}

Two weighing matrices, $H$ and $K$, are {\it equivalent} if $H = PKQ$, where $P$ and $Q$ are unimodular permutation matrices (i.e., each row/column has exactly one nonzero unimodular entry). We use the notation $H \cong K$.

\begin{definition}
 Let $W$ be a weighing matrix of order $n$ and weight $w$. If $W \cong W_1 \oplus W_2$ for some $W_1$ and $W_2$ of order strictly less than $n$, then $W$ is said to be {\bf decomposable}. We may write $W$ in such a way that $W=W_1 \oplus W_2 \oplus \cdots \oplus W_k$ where each $W_i$ is not decomposable of order $n_i$ such that $n_j \leq n_{j+1}$ for $1 \leq j < k$. The {\bf block structure} of $W$ is the matrix $$J_{n_1} \oplus J_{n_2} \oplus \cdots \oplus J_{n_k},$$ where $J_n$ is the all ones matrix of order $n$.
\end{definition}

Determining if two weighing matrices are equivalent is a relatively complex problem, and as of today, there are no efficient algorithms to determine equivalence. Determining if two weighing matrices have the same block structure, however, is a much simpler problem as we see in the next lemma.

\begin{lemma}\label{lem:block-complexity}
 The block structure of a weighing matrix can be determined in $O(n^3)$.
 \begin{proof}
  Given a weighing matrix $W$ of order $n$, let $G$ be the graph on $n$ vertices with an edge between $i$ and $j$ if and only if at least one nonzero entry in row $i$ is in the same column as a nonzero entry in row $j$. Two rows of $W$ are in the same non-decomposable block if and only if there is a path between the corresponding nodes in $G$. Thus, a non-decomposable block of $W$ can be found by taking the rows corresponding to all vertices in any connected component of $G$ and removing any columns that only have zeroes. The number of non-decomposable blocks of $W$ is the number of connected components of $G$. By placing the number of vertices in each non-decomposable block into a list and sorting that list (say we now have $n_1,n_2,\cdots,n_k$), we have that the block structure of $W$ is $$J_{n_1} \oplus J_{n_2} \oplus \cdots \oplus J_{n_k}.$$

 This process has three steps: First, we must build the graph. This can be done in $O(n^3)$ by looking at all pairs of rows and examining each column. Then, we determine the number of connected components, which takes $O(n^2)$ via depth first search. Finally, we sort the list in $O(n\log n)$ for a total complexity of $O(n^3)$.
 \end{proof}
\end{lemma}

 It is noteworthy to point out that the asymptotic bound in Lemma \ref{lem:block-complexity} is not tight. When constructing $G$ in the proof of Lemma \ref{lem:block-complexity} can be done by multiplying $|W|$ by $|W|^T$, where $|W| = [|w_{ij}|]$. The nonzero entries in $|W||W|^T$ signify an edge in $G$. As of today, matrix multiplication can be done in $O(n^{2.3727})$, but in general, due to the fact that we are only concerned with the fact that an entry is nonzero, we can apply bit operations to make the $O(n^3)$ algorithm significantly faster in practice.

\begin{proposition} \label{prop:blocks}
 If two weighing matrices (say $H$ and $K$) of the same weight have the same block structure, then $H$ is unbiased with $K$ if and only if each non-decomposable block of $H$ is unbiased with the corresponding non-decomposable block of $K$.
 \begin{proof}
  This is easily seen by noting that $$(H_1 \oplus \cdots \oplus H_m)(K_1 \oplus \cdots \oplus K_m)^* = (H_1K_1^* \oplus \cdots \oplus H_mK_m^*).$$
 \end{proof}

\end{proposition}

\begin{proposition} \label{prop:blocks_upper}
 If every matrix in a set of mutually unbiased weighing matrices has the same block structure, then that set's size is restricted by each individual non-decomposable block's upper bound.
 \begin{proof} This follows from Proposition \ref{prop:blocks}.\end{proof}
\end{proposition}

The following two theorems from Calderbank et al. \cite{calderbank97} are very important results that we will be using.

\begin{theorem}(\cite[Equation 5.9]{calderbank97})\label{th:vec_upper}
  Let $V \subset \mathbb{C}^n$ be a set of unit vectors. If $|\langle v,w\rangle| \in \left\{0,\alpha\right\}$ for all $v,w\in V$, $v\neq w$, where $\alpha \in \mathbb{R}$ and $0 < \alpha < 1$, then

\begin{equation}|V|\leq n\binom{n+1}{2}.\end{equation} Moreover, \begin{equation}|V|\leq \frac{n(n+1)(1-\alpha^2)}{2-(n+1)\alpha^2}\end{equation} if the denominator is positive.
\end{theorem}

\begin{theorem}(\cite[Equation 3.7  and 3.9]{calderbank97})\label{th:vec_upper_real}
  If all of the entries of $V$ in Theorem \ref{th:vec_upper} are real, then \begin{equation}|V|\leq \binom{n+2}{3}.\end{equation} Moreover, \begin{equation}|V|\leq \frac{n(n+2)(1-\alpha^2)}{3-(n+2)\alpha^2}\end{equation} if the denominator is positive.
\end{theorem}

It is important to note that in most cases, the second upper bound given in each theorem is smaller than the first, but not always. For example, if we are looking for real vectors with $n=9$ and $\alpha = \frac12$, the first bound gives us $|V| \leq 165$ whereas the second bound gives us $|V| \leq 297$.

The following are immediate corollaries to the previous two theorems.

\begin{corollary}
\label{cor:mat_upper}
 Let $\mathcal{W}=\{W_1,\cdots,W_m\}$ be a set of mutually unbiased weighing matrices of order $n$ and weight $w$. Then we have that the size of $\mathcal{W}$ is bounded above by

  \begin{equation}\label{eq:aub} \frac{(n-1)(n+2)}{2}.\end{equation}
  Moreover, if $2w-(n+1) > 0$, then it is bounded above by
  \begin{equation}\label{eq:tub} \frac{w(n-1)}{2w-(n+1)}.\end{equation}
  
  \begin{proof}
    Define $V$ to be the set of all rows of $\frac{1}{\sqrt{w}}W_1, \cdots, \frac{1}{\sqrt{w}}W_m$ (note that $|V| = mn$). Since $\mathcal{W}$ is a set of mutually unbiased weighing matrices, we may set $\alpha=\frac{1}{\sqrt{w}}$. Moreover, note that since all vectors in $V$ come from a weighing matrix of weight $w$, we may add the rows of the identity matrix to $V$ without disrupting the bi-angularity. By applying Theorem \ref{th:vec_upper} to $V$ (with the identity rows), we arrive at the desired results.
   \end{proof}
\end{corollary}

\begin{corollary}
\label{cor:mat_upper_real}
 Let $\mathcal{W}=\{W_1,\cdots,W_m\}$ be a set of real mutually unbiased weighing matrices of order $n$ and weight $w$. Then we have that the size of $\mathcal{W}$ is bounded above by

  \begin{equation}\label{eq:aub_real}\frac{(n-1)(n+4)}{6}.\end{equation}
  Moreover, if $3w-(n+2) > 0$, then it is bounded above by
  \begin{equation}\label{eq:tub_real} \frac{w(n-1)}{3w-(n+2)}.\end{equation}

  \begin{proof}
   Similar to Corollary \ref{cor:mat_upper}.
  \end{proof}

\end{corollary}

\section{Mutually Unbiased Weighing Matrices} \label{sec:main}

\subsection{Computer Search}

With unit weighing matrices, an exhaustive computer search is impractical, if not impossible, to perform since each nonzero entry in each matrix has infinitely many choices. To this end, we restricted the entries to small roots of unity in our computer searches. For each type of matrix, we searched for matrices over the $m^{th}$ roots of unity, with $m \leq 24$. As one observes from Table \ref{table:bounds}, the $12^{th}$ roots of unity seem to be the largest group needed to find some maximal sets. Many of the maximal sets that we found do not match the upper bound given in Corollary \ref{cor:mat_upper}. For many cases, we prove smaller upper bounds.

\begin{table}[H]
\caption{We compare the theoretic upper bound given in Corollary \ref{cor:mat_upper} to the results of both our computer searches and any improved (i.e., smaller) upper bounds we have found. The highlighted rows signify cases where the smallest upper bound and largest lower bound do not meet. Note that $UW(6,6)$ is the most highly sought after set of matrices.}
\centering
\begin{tabular}{@{}cccllccc@{}}
 \toprule\label{table:bounds}
 Type & & \multicolumn{2}{@{}c@{}}{Upper Bounds} & & \multicolumn{2}{@{}c@{}}{Examples Found} \\
\cmidrule{2-4} \cmidrule{6-7}
 & & Corollary \ref{cor:mat_upper} & \multicolumn{1}{@{}c@{}}{Smallest} & & Largest Set & Root of Unity \\
 \midrule
 UW(2,2)       && 2              & 2                          && 2  &  4  \\
 UW(3,2)       && 5              & 0 (See \cite{unit-weigh13})   && 0  & --  \\
 UW(3,3)       && 3              & 3                          && 3  &  3  \\
 UW(4,2)       && 9              & 2 (Lemma \ref{lem:even-2}) && 2  &  4  \\
 UW(4,3)       && 9              & 9                          && 9  &  6  \\
 UW(4,4)       && 4              & 4                          && 4  &  4  \\
 UW(5,2)       && 14             & 0 (See \cite{unit-weigh13})   && 0  & -- \\
 UW(5,3)       && 14             & 0 (See \cite{unit-weigh13})   && 0  & -- \\
 UW(5,4)       && 8              & 5 (Theorem \ref{th:cw_5_4})&& 5  &  6  \\
 UW(5,5)       && 5              & 5                          && 5  &  5  \\
 UW(6,2)       && 20             & 2 (Lemma \ref{lem:even-2}) && 2  &  4  \\
 UW(6,3)       && 20             & 3 (Theorem \ref{th:cw_n_3})&& 3  &  3  \\
 UW(6,4)       && 20             & 20                         && 20 &  6  \\
 \rowcolor{Gray}
 UW(6,5)       && $\frac{25}{3}$ & 8                          && 2  & 12  \\
 \rowcolor{Gray}
 UW(6,6)       && 6              & 6                          && 2  & 12  \\
 UW(7,2)       && 27             & 0 (See \cite{unit-weigh13})   && 0  & -- \\
 UW(7,3)       && 27             & 3 (Theorem \ref{th:cw_n_3})&& 3  &  6  \\
 UW(7,4)       && 27             & 8 (Theorem \ref{th:cw_7_4})&& 8  &  2  \\
 UW(7,5)       && 15             & 0 (See \cite{unit-weigh13})   && 0  & -- \\
 \rowcolor{Gray}
 UW(7,6)       && 9              & 9                          && 0  &  --   \\
 UW(7,7)       && 7              & 7                          && 7  &  7  \\
 \bottomrule
\end{tabular}
\end{table}

Mutually unbiased unit Hadamard matrices have been extensively studied for prime power orders. A proof of the following Theorem can be found in \cite{Bengtsson_2007}.

\begin{theorem}
\label{th:pr-pow}
 For any prime power $p$, there exists a full set of $p$ mutually unbiased (Butson) Hadamard matrices $UW(p,p)$.
\end{theorem}

\subsection{Upper bound for Mutually Unbiased Weighing Matrices of Weight 2}
In \cite[Theorem 10]{unit-weigh13}, we proved that $UW(n,2)$ do not exist for odd orders. For $n$ even, we have the following.

\begin{lemma} \label{lem:even-2}
 Let $n$ be even. Then there are at most 2 mutually unbiased weighing matrices of order $n$ and weight 2.
 \begin{proof}
  Say we have a set of mutually unbiased weighing matrices of the appropriate order and weight. From \cite{unit-weigh13}, we know that one of the matrices may be transformed into $$\left(
\begin{array}{cc}
 1 & 1 \\
 1 & -\\
\end{array}
\right) \otimes I_{n/_2}
.$$ Permute the rows of the second matrix so that there is a nonzero in the top-left entry. The second entry in the top row must be nonzero, otherwise the inner product of the top row of the first and second matrices will be neither 0 nor $\sqrt{2}$. Continue this argument so that the block structure is the same between all matrices in the set of unbiased weighing matrices. By applying Corollary \ref{cor:mat_upper} to Proposition \ref{prop:blocks_upper}, we have our result.
 \end{proof}

\end{lemma}

\subsection{Upper bound for Mutually Unbiased Weighing Matrices of Weight 3}
\begin{lemma} \label{lem:cw_n_3_blocks}
 A $UW(n,3)$, $H$, is unbiased with $K$ if and only if $K$ has the same block structure as $H$.

 \begin{proof}
  From \cite[Theorem 12]{unit-weigh13}, we know that $H$ may be transformed into a matrix of the following form:

$$\left(
\begin{array}{ccc}
 1 & 1       &1 \\
 1 & a       &\overline{a} \\
 1 & \overline{a} &a       
\end{array}
\right) \oplus \cdots \oplus \left(
\begin{array}{ccc}
 1 & 1       &1 \\
 1 & a       &\overline{a} \\
 1 & \overline{a} &a       
\end{array}
\right) \oplus 
\left(
\begin{array}{cccc}
 1 & 1 &1 &0 \\
 1 & - &0 &1 \\
 1 & 0 &- &- \\
 0 & 1 &- &1
\end{array}
\right)
\oplus \cdots \oplus
\left(
\begin{array}{cccc}
 1 & 1 &1 &0 \\
 1 & - &0 &1 \\
 1 & 0 &- &- \\
 0 & 1 &- &1
\end{array}
\right).
$$

We may assume that the first 3 rows of $K$ have a 1 in the first column by appropriate row permutations.

Assume that the top left block in $H$ is a $UW(3,3)$. If columns 2 and 3 of $K$ are both zero in any of the first 3 rows, then the inner product of row 1 in $H$ and that row will give us a unimodular number, not having absolute value 0 or $\sqrt3$. If exactly one of the entries in columns 2 and 3 are nonzero, then there must be a third nonzero in one of the last $n-3$ columns. Taking the inner product of this row and an appropriate row in $H$, there is another unimodular number, causing the same contradiction as above. Thus,  in these first three rows of $K$, each must have exactly 3 nonzero entries in the first three columns (ie. a $UW(3,3)$).

Now assume that the top left block in $H$ is a $UW(4,3)$. If columns 2,3 and 4 are all zero in any of the first 3 rows, then the inner product of row 1 in $H$ and that row will give us a unimodular number. If there is exactly 1 nonzero in columns 2,3 and 4, then the inner product of that row and the fourth row of $H$ will be unimodular. Thus, we know that in the first 3 rows of $K$, all 3 nonzero entries must appear in the first four columns.

We will now show that the first zero in these rows will not be in the same column. Assume that one column has at least 2 zeroes. This means that at least one of columns 2,3 and 4 will be complete (i.e., no more nonzero entries may go into that column). Column 1 is already complete, so in our fourth row, there is either 1 or 2 nonzeroes in the first 3 columns. By taking the inner product of the fourth row of $K$ by the appropriate row in $H$, we will get a unimodular number. Thus, the first zero in the first 4 rows must be in different columns (note that the first zero in row 4 must be in column 1). Furthermore, through appropriate row permutations and negations, the second entry in row 4 must be a 1. The next two entries are clearly nonzero or there is 1-orthogonality within $K$. Thus, in the first 4 rows of $K$, the three nonzero entries must appear in the first 4 rows, with the first zeroes of the rows in different columns (i.e., a $UW(4,3)$).

Once we know that the top left block of $H$ and $K$ are the same, if we examine the bottom right $(n-3) \times (n-3)$ or $(n-4) \times (n-4)$ block, we have a $UW(n-3,3)$ or $UW(n-4,3)$, and we can recursively use the same argument to obtain the desired result.
\end{proof}

\end{lemma}

\begin{theorem} \label{th:cw_n_3}
 The upper bound on the number of MUWM of the form $UW(n,3)$ is:
  $$\begin{cases}
   3 & \text{if } n \not\equiv 0 \mod 4 \\
   9 & \text{if } n \equiv 0 \mod 4
  \end{cases}$$
\\ where $n \in \{3,4\} \cup \{k:k\geq 6\}$.
\begin{proof}

Using Lemma \ref{lem:cw_n_3_blocks} with Proposition \ref{prop:blocks_upper} and the fact that the upper bound for $UW(3,3)$ is 3 and $UW(4,3)$ is 9 via Corollary \ref{cor:mat_upper}, we have that if the matrix contains a $UW(3,3)$ in its block structure, then it acts as a limiting factor, causing the upper bound to be 3. Otherwise, it is 9, which can only occur when $n$ is a multiple of 4.
\end{proof}
\end{theorem}

\begin{corollary} \label{cor:cw_n_3}
 The upper bound given in Theorem \ref{th:cw_n_3} is tight for all $n \in \{3,4\} \cup \{k:k\geq 6\}$.
\begin{proof}
 A computer search has shown the bounds to be tight for $UW(4,3)$ and the bound for $UW(3,3)$ is attained through Theorem \ref{th:pr-pow}. We may construct the $UW(n,3)$ by adjoining the appropriate amount of $UW(4,3)$ and $UW(3,3)$ together along the main diagonals. If $n$ is a multiple of 4, use only $UW(4,3)$s along the main diagonal. Otherwise, it does not matter which blocks are used. A simple induction will show that every integer larger than 5 may be written in the form of $3m+4l$.
\end{proof}
\end{corollary}

\subsection{Upper bound for Mutually Unbiased Weighing Matrices of Weight 4}

\subsubsection{UW(5,4)}

\begin{lemma}\label{lem:cw_5_4}
 Let $W$ be a unit weighing matrix that is unbiased with $W_5$, then every nonzero entry in $W$ is a sixth root of unity. $W_5$ is given as follows:
 Let $W$ be a unit weighing matrix that is unbiased with the following matrix:
$$
W_5 = \left(
\begin{array}{ccccc}
 1 & 1            &1            &1            &0 \\
 1 & a            &\overline{a} &0            &1 \\
 1 & \overline{a} &0            &a            &\overline{a} \\
 1 & 0            &a            &\overline{a} &a \\
 0 & 1            &\overline{a} &a            &a
\end{array}
\right)
$$
where $a = e^{i\frac{2\pi}{3}}$.

 \begin{proof}
  Since $W_5W^* = 2L$ for some weighing matrix $L$, we know that each row of $W$ must be orthogonal with exactly one row of $W_5$. Moreover, we may permute the rows of $W$ so that row $i$ is orthogonal with row $i$ of $W_5$. We know that the first nonzero entry in each row of $W$ may be a one. Using the definition of $m$-orthogonality and the results given in \cite[Section 3]{unit-weigh13}, we can determine that there are at most 11 {\it different} rows that are orthogonal to each of the rows of $W_5$, each with exactly one free variable.

  Let $b$ be an arbitrary unimodular number and $\alpha$ a primitive third root of unity. The four main observations that are used in each line of the proof are:
\begin{enumerate}
 \item[(O1)] $|1-\alpha+b|=2 \implies b=\pm\overline{\alpha}$,
 \item[(O2)] $|1+\alpha+b|=2 \implies b=-\overline{\alpha}$,
 \item[(O3)] $|3+b|=2 \implies b=-1$,
 \item[(O4)] $1+\alpha+\overline{\alpha}=0$.
\end{enumerate}

  We will examine all candidates for row 1 of $W$. There are only 11 different candidates (up to a free variable), they are:
     $$
     \begin{array}{rccccc}
	(A)  & 1 & -       & b       & -b      & 0 \\
	(B)  & 1 & b       & -       & -b      & 0 \\
	(C)  & 1 & b       & -b      & -       & 0 \\
	(D)  & 1 & a       & \overline{a} & 0       & b \\
	(E)  & 1 & \overline{a} & a       & 0       & b \\
	(F)  & 1 & a       & 0       & \overline{a} & b \\
	(G)  & 1 & \overline{a} & 0       & a       & b \\
	(H)  & 1 & 0       & a       & \overline{a} & b \\
	(I)  & 1 & 0       & \overline{a} & a       & b \\
	(J)  & 0 & 1       & a       & \overline{a} & b \\
	(K)  & 0 & 1       & \overline{a} & a       & b \\
     \end{array}
     $$

     For each candidate, we will show that in order to be unbiased with the other four rows of $W_5$, the free variable must be a sixth root of unity. In some cases, we will show that the row cannot be unbiased with a specific row of $W_5$. To avoid a lengthy proof, we only give three examples.

     \begin{enumerate}
       \item[(A)] By taking the complex inner product with row 2 of $W_5$, we have that $|1-a+\overline{ab}|=2$. By using $(O1)$, we have that $\overline{ab} = \pm\overline{a}$ which implies  that $b = \pm1$. Thus, all entries in the candidate row are sixth roots of unity.
       \item[(G)] By taking the complex inner product with row 3 of $W_5$, we have that $|1+1+1+\overline{ab}|=2$. By using $(O3)$, we have that $\overline{ab} = -1$ which implies that $b = -\overline{a}$. Thus, all entries in the candidate row are sixth roots of unity.
       \item[(J)] By taking the complex inner product with row 5 of $W_5$, we have that $|1+a+\overline{a}+a\overline{b}|=2$. By using $(O4)$, we have that $|a\overline{b}| = 2$ which implies that $|b| = 2$, which is a contradiction since $b$ is a unimodular number. Thus, $(J)$ cannot be unbiased with row 5, so it may not be the row that is orthogonal with row $1$ of $W_5$.
     \end{enumerate}

      For each of the 5 rows of $W_5$, there are 11 {\it different} candidates for each row (each with exactly one free variable). In each case, the free variable is shown to be a sixth root of unity or have absolute value 2 (as in the examples above).
 \end{proof}
\end{lemma}

\begin{theorem} \label{th:cw_5_4}
 The largest number of mutally unbiased weighing matrices of the form $UW(5,4)$ is 5.
 \begin{proof}
  In \cite[Lemma 15]{unit-weigh13}, it is proven that all $UW(5,4)$ are equivalent to $W_5$ given in Lemma \ref{lem:cw_5_4}. Thus, given a set of mutually unbiased weighing matrices, we may permute and negate the rows and columns of the matrices in such a way that one of them is $W_5$. By Lemma \ref{lem:cw_5_4}, we know that any matrix that is unbiased with $W_5$ must only contain 0 and the sixth roots of unity. An exhaustive computer search was done over these entries, which reveiled that the maximal set using only the sixth root of unity contains 5 elements. These matrices are included in Appendix \ref{app:sets}.
 \end{proof}
\end{theorem}

\subsubsection{UW(7,4)}

\begin{lemma} \label{lem:cw_7_4}
 Let $W$ be a unit weighing matrix that is unbiased with $W_7$, then every nonzero entry in $W$ is real. $W_7$ is given as follows:
$$W_7=\left(\begin{array}{c}
\Zp\Zp\Zp\Zp\Zz\Zz\Zz\\
\Zp\Zm\Zz\Zz\Zp\Zp\Zz\\
\Zp\Zz\Zm\Zz\Zm\Zz\Zp\\
\Zp\Zz\Zz\Zm\Zz\Zm\Zm\\
\Zz\Zp\Zm\Zz\Zz\Zp\Zm\\
\Zz\Zp\Zz\Zm\Zp\Zz\Zp\\
\Zz\Zz\Zp\Zm\Zm\Zp\Zz
\end{array}\right).$$

 \begin{proof}
  We can easily see that there are only $\binom73=35$ possible zero placements that are valid in a row of $W$. Similar to the proof of Lemma \ref{lem:cw_5_4}, we will only show a couple cases, as the rest follow similarly. Let $a,b,c$ be arbitrary unimodular numbers.

  \begin{enumerate}
   \item[(A)] $\left(\begin{array}{ccccccc}1 & a & b & c & 0 & 0 & 0\end{array}\right)$
     \begin{itemize}
      \item Taking the complex inner product with row 2 of $W_7$, we have that $|1+a| \in \{0,2\}$ which implies $a \in \{\pm1\}$.
      \item Taking the complex inner product with row 3 of $W_7$, we have that $|1+b| \in \{0,2\}$ which implies $b \in \{\pm1\}$.
      \item Taking the complex inner product with row 4 of $W_7$, we have that $|1+c| \in \{0,2\}$ which implies $c \in \{\pm1\}$.
     \end{itemize}
   \item[(B)] $\left(\begin{array}{ccccccc}1 & a & b & 0 & c & 0 & 0\end{array}\right)$
     \begin{itemize}
      \item Taking the complex inner product with row 4 of $W_7$, we have that $|1| \in \{0,2\}$ which is clearly a contradiction.
     \end{itemize}
  \end{enumerate}
 \end{proof}
\end{lemma}

Of particular note, the only rows that do not cause a contradiction are those 7 rows which have the same zero placement as $W_7$.

\begin{theorem} \label{th:cw_7_4}
 The maximum number of mutually unbiased weighing matrices of order 7 and weight 4 is 8.
 \begin{proof}
  Similar to the proof of Theorem \ref{th:cw_5_4}, one matrix in the set may be transformed into the real weighing matrix $W_7$ given Lemma \ref{lem:cw_7_4}. Every $UW(7,4)$ is equivalent to this matrix (see \cite[Section 3.4]{unit-weigh13}). By Lemma \ref{lem:cw_7_4}, every weighing matrix equivalent to $W_7$ must also be real, so we may use Corollary \ref{cor:mat_upper_real} to provide us with this bound.
 \end{proof}
\end{theorem}

Using a computer search, we find the eight real mutually unbiased weighing matrices $W(7,4)$ given in Appendix \ref{app:sets}. This achieves the real upper bound given by Corollary \ref{cor:mat_upper_real}.
By Theorem \ref{th:cw_7_4}, this is also the maximal set of $UW(7,4)$, despite not achieving the upper bound of 24 given by Corollary \ref{cor:mat_upper}.

\subsubsection{UW(8,4)}
\begin{theorem}\label{th:w_8_4}
 The maximum number of real mutually unbiased weighing matrices of order 8 and weight 4 is 14.
 \begin{proof}
  A set of size 14 $W(8,4)$ has been generated in Appendix \ref{app:sets}. This meets the upper bound given by Corollary \ref{cor:mat_upper_real}.
 \end{proof}
\end{theorem}

Further inverstigations into $UW(8,4)$ using large roots of unity have proven fruitless. Odd roots of unity produce maximal sets smaller than that of the real case, and even roots of unity become computationally infeasible after the fourth root of unity, which returns the set of $W(8,4)$ as the maximal set of mutually unbiased weighing matrices.

\section{Unbiased Hadamard Matrices}\label{sec:weakly-unbiased}

Thus far, we have only examined a very special case of unbiasedness. Our selection of the values of $\alpha$ in (\ref{eq:tub}) and (\ref{eq:tub_real}) make it possible to  append the identity to the set of weighing matrices. More preciesly, considering each row of all weighing matrices in a set of mutually unbiased weighing matrices of order $n$ and the rows of the identity matrix of order $n$  as vectors in $\mathbb{R}^n$ or $\mathbb{C}^n$, they form a class of bi-angular vectors. We now make a different selction for the value of $\alpha$ in (\ref{eq:tub_real}) in such a way that it is no longer possible to add the identity matrix and preserve the bi-angularity. Below, we give an example of a set of eight Hadamard matrices of order 8 that form a bi-angular set of vectors in $\mathbb{R}^8$, but no rows of the identity matrix can be added to the set and preserve bi-angularity. In the following set, $\alpha=\frac12$, but if the identity is added, it would introduce the inner product of $\frac{1}{\sqrt{8}}$ (up to absolute value) and the bi-angularity of the lines disappear.

\begin{table}[H]\centering\caption{8 mutually unbiased Hadamard matrices with $\alpha=\frac12$}
\begin{tabular}{@{}cccc@{}}\toprule\label{table:H8}
$\left(\begin{array}{c}
 \Zp\Zp\Zp\Zp\Zp\Zp\Zp\Zp \\
 \Zp\Zp\Zm\Zp\Zm\Zm\Zp\Zm \\
 \Zp\Zm\Zm\Zp\Zp\Zm\Zm\Zp \\
 \Zp\Zm\Zm\Zm\Zm\Zp\Zp\Zp \\
 \Zp\Zp\Zm\Zm\Zp\Zp\Zm\Zm \\
 \Zp\Zp\Zp\Zm\Zm\Zm\Zm\Zp \\
 \Zp\Zm\Zp\Zp\Zm\Zp\Zm\Zm \\
 \Zp\Zm\Zp\Zm\Zp\Zm\Zp\Zm
\end{array}\right)$ & $
\left(\begin{array}{c}
 \Zp\Zp\Zp\Zm\Zp\Zm\Zp\Zp \\
 \Zp\Zm\Zp\Zp\Zp\Zp\Zp\Zm \\
 \Zp\Zm\Zm\Zp\Zm\Zm\Zp\Zp \\
 \Zp\Zp\Zm\Zp\Zp\Zm\Zm\Zm \\
 \Zp\Zp\Zm\Zm\Zm\Zp\Zp\Zm \\
 \Zp\Zm\Zp\Zm\Zm\Zm\Zm\Zm \\
 \Zp\Zm\Zm\Zm\Zp\Zp\Zm\Zp \\
 \Zp\Zp\Zp\Zp\Zm\Zp\Zm\Zp
\end{array}\right)$ & $
\left(\begin{array}{c}
 \Zp\Zp\Zm\Zm\Zm\Zp\Zm\Zp \\
 \Zp\Zm\Zm\Zm\Zp\Zp\Zp\Zm \\
 \Zp\Zm\Zp\Zm\Zm\Zm\Zp\Zp \\
 \Zp\Zp\Zp\Zp\Zm\Zp\Zp\Zm \\
 \Zp\Zp\Zp\Zm\Zp\Zm\Zm\Zm \\
 \Zp\Zm\Zp\Zp\Zp\Zp\Zm\Zp \\
 \Zp\Zm\Zm\Zp\Zm\Zm\Zm\Zm \\
 \Zp\Zp\Zm\Zp\Zp\Zm\Zp\Zp
\end{array}\right)$ & $
\left(\begin{array}{c}
 \Zp\Zm\Zm\Zm\Zm\Zp\Zm\Zm \\
 \Zp\Zp\Zp\Zm\Zm\Zm\Zp\Zm \\
 \Zp\Zp\Zm\Zm\Zp\Zp\Zp\Zp \\
 \Zp\Zm\Zp\Zp\Zm\Zp\Zp\Zp \\
 \Zp\Zm\Zp\Zm\Zp\Zm\Zm\Zp \\
 \Zp\Zp\Zp\Zp\Zp\Zp\Zm\Zm \\
 \Zp\Zp\Zm\Zp\Zm\Zm\Zm\Zp \\
 \Zp\Zm\Zm\Zp\Zp\Zm\Zp\Zm
\end{array}\right)$ \\[30pt] $
\left(\begin{array}{c}
 \Zp\Zm\Zp\Zm\Zm\Zp\Zm\Zp \\
 \Zp\Zp\Zp\Zm\Zp\Zp\Zp\Zm \\
 \Zp\Zp\Zm\Zm\Zm\Zm\Zp\Zp \\
 \Zp\Zm\Zm\Zp\Zm\Zp\Zp\Zm \\
 \Zp\Zm\Zm\Zm\Zp\Zm\Zm\Zm \\
 \Zp\Zp\Zp\Zp\Zm\Zm\Zm\Zm \\
 \Zp\Zp\Zm\Zp\Zp\Zp\Zm\Zp \\
 \Zp\Zm\Zp\Zp\Zp\Zm\Zp\Zp
\end{array}\right)$ & $
\left(\begin{array}{c}
 \Zp\Zm\Zm\Zp\Zm\Zp\Zm\Zp \\
 \Zp\Zm\Zm\Zm\Zp\Zm\Zp\Zp \\
 \Zp\Zp\Zm\Zp\Zp\Zp\Zp\Zm \\
 \Zp\Zp\Zp\Zp\Zm\Zm\Zp\Zp \\
 \Zp\Zm\Zp\Zp\Zp\Zm\Zm\Zm \\
 \Zp\Zm\Zp\Zm\Zm\Zp\Zp\Zm \\
 \Zp\Zp\Zp\Zm\Zp\Zp\Zm\Zp \\
 \Zp\Zp\Zm\Zm\Zm\Zm\Zm\Zm
\end{array}\right)$ & $
\left(\begin{array}{c}
 \Zp\Zp\Zp\Zm\Zm\Zp\Zm\Zm \\
 \Zp\Zm\Zp\Zm\Zp\Zp\Zp\Zp \\
 \Zp\Zm\Zm\Zm\Zm\Zm\Zp\Zm \\
 \Zp\Zp\Zm\Zp\Zm\Zp\Zp\Zp \\
 \Zp\Zp\Zm\Zm\Zp\Zm\Zm\Zp \\
 \Zp\Zm\Zm\Zp\Zp\Zp\Zm\Zm \\
 \Zp\Zm\Zp\Zp\Zm\Zm\Zm\Zp \\
 \Zp\Zp\Zp\Zp\Zp\Zm\Zp\Zm
\end{array}\right)$ & $
\left(\begin{array}{c}
 \Zp\Zp\Zm\Zm\Zp\Zm\Zp\Zm \\
 \Zp\Zm\Zp\Zp\Zm\Zm\Zp\Zm \\
 \Zp\Zm\Zm\Zp\Zp\Zp\Zp\Zp \\
 \Zp\Zp\Zp\Zp\Zp\Zm\Zm\Zp \\
 \Zp\Zp\Zp\Zm\Zm\Zp\Zp\Zp \\
 \Zp\Zm\Zp\Zm\Zp\Zp\Zm\Zm \\
 \Zp\Zm\Zm\Zm\Zm\Zm\Zm\Zp \\
 \Zp\Zp\Zm\Zp\Zm\Zp\Zm\Zm
\end{array}\right)$\\
\bottomrule
\end{tabular}
\end{table}

The rows of these matrices are generated from the BCH-code of length 7 with weight distribution $\left\{(0,1),(2,21),(4,
35),(6,7)\right\}$ (see \cite{BCH-BC,BCH-H} for more information about BCH-codes). Once the codewords are generated, we append a column of zeroes, then perform the following operation onto each entry of the codewords:

$$f(i) = \begin{cases}
          1 & \text{ if } i = 0 \\
         -1 & \text{ if } i = 1
         \end{cases}.$$
         
We were also able to generate 32 Hadamard matrices of order 32 which have inner products in $\left\{0,\pm8 \right\}$ through a similar process. The weight distribution of the order 32 matrices is $\left\{(0,1), (12, 310), (16, 527), (20, 186)\right\}$. The partition of the vectors into Hadamard matrices is shown in Appendix \ref{app:H32}.

In an attempt to continue this, we have generated the $128^2$ codewords from the BCH-code of order 127, but were not able to partition them into the 128 Hadamard matrices needed due to computer memory restrictions. The inner products between the vectors are all in $\left\{0,\pm16\right\}$. We do believe that this set of vectors contains the needed ingredients to make the Hadamard matrices required. Moreover, we pose the following

\begin{conjecture}\label{conj:real-hada}
 Let $n = 2^{2k+1}$. Then there exists a set of $n$ real Hadamard matrices, $\left\{H_1, H_2, \cdots , H_n\right\}$, so that the entries of $H_iH_j^t$ ($i \neq j$) contain exactly two elements, $0$ and $2^{k+1}$ (up to absolute value).
\end{conjecture}

It is important to note that the number of vectors found through Conjecture \ref{conj:real-hada} is usually less than the bound given in Theorem \ref{th:vec_upper_real}. We believe that the upper bound is too high in this case because the vectors are {\it flat} (i.e., all contain entries that have the same absolute value). In fact, it is our belief that when the restriction of {\it flatness} is imposed on a set of vectors or matrices, a much smaller general upper bound should be possible.

Using the terminology from \cite{much10}, these matrices form a set of {\it weakly unbiased Hadamard matrices}. However, it is important to note that the matrices formed here are a very special kind of unbiased Hadamard matrices since the entire set of vectors forms a set of bi-angular lines (whereas the vectors from \cite{much10} give tri-angular lines). These matrices seem to form very nice combinatorial objects, which are discussed in further detail in \cite{biangular_lines}.

{\bf Acknowledgments:} The authors wish to extend their gratitude to Professor Masaaki Harada for his help in locating the codes used in section \ref{sec:weakly-unbiased}. The authors also wish to thank Professor Kevin Grant for allowing the use of his NSERC funded computer, {\it hera}, for many of the computations carried out in this article.

\appendix

\section{Sets Attaining the Smallest Upper Bound}\label{app:sets}

This section includes a library of sets of weighing matrices whose size equal the smallest upper bound that is known. To save space, we define $a := e^{2\pi i / 3}$ and $b := e^{2\pi i / 6}$.

\begin{table}[H]\centering\caption{9 mutually unbiased weighing matrices of order 4 and weight 3, $UW(4,3)$.}
\begin{tabular}{@{}ccccc@{}}\toprule\label{table:UW4_3}
$
\left(\begin{array}{c}
\Zp\Zp\Zp\Zz\\
\Zp\Zm\Zz\Zp\\
\Zp\Zz\Zm\Zm\\
\Zz\Zp\Zm\Zp\\
\end{array}\right)$ & $
\left(\begin{array}{c}
\Zp\Zp\Za\Zz\\
\Zp\Zm\Zz\Za\\
\Zp\Zz\Zb\Zb\\
\Zz\Zp\Zb\Za\\
\end{array}\right)$ & $
\left(\begin{array}{c}
\Zp\Zp\ZA\Zz\\
\Zp\Zm\Zz\ZA\\
\Zp\Zz\ZB\ZB\\
\Zz\Zp\ZB\ZA\\
\end{array}\right)$ & $
\left(\begin{array}{c}
\Zp\ZB\Zz\Zp\\
\Zp\ZA\Za\Zz\\
\Zp\Zz\Zb\Zm\\
\Zz\Zp\ZB\Za\\
\end{array}\right)$ & $
\left(\begin{array}{c}
\Zp\ZB\Zz\Za\\
\Zp\ZA\ZA\Zz\\
\Zp\Zz\ZB\Zb\\
\Zz\Zp\Zm\ZA\\
\end{array}\right)$ \\ $
\left(\begin{array}{c}
\Zp\ZB\Zz\ZA\\
\Zp\ZA\Zp\Zz\\
\Zp\Zz\Zm\ZB\\
\Zz\Zp\Zb\Zp\\
\end{array}\right)$ & $
\left(\begin{array}{c}
\Zp\Za\Zp\Zz\\
\Zp\Zb\Zz\Za\\
\Zp\Zz\Zm\Zb\\
\Zz\Zp\ZB\Zp\\
\end{array}\right)$ & $
\left(\begin{array}{c}
\Zp\Za\Za\Zz\\
\Zp\Zb\Zz\ZA\\
\Zp\Zz\Zb\ZB\\
\Zz\Zp\Zm\Za\\
\end{array}\right)$ & $
\left(\begin{array}{c}
\Zp\Za\ZA\Zz\\
\Zp\Zb\Zz\Zp\\
\Zp\Zz\ZB\Zm\\
\Zz\Zp\Zb\ZA\\
\end{array}\right)
$\\
\bottomrule
\end{tabular}
\end{table}

\begin{table}[H]\centering\caption{5 mutually unbiased weighing matrices of order 5 and weight 4, $UW(5,4)$.}
\begin{tabular}{@{}ccccc@{}}\toprule\label{table:UW5_4}
$
\left(\begin{array}{c}
\Zp\Zp\Zp\Zp\Zz\\
\Zp\Za\ZA\Zz\Zp\\
\Zp\ZA\Zz\Za\ZA\\
\Zp\Zz\Za\ZA\Za\\
\Zz\Zp\ZA\Za\Za\\
\end{array}\right)$ & $
\left(\begin{array}{c}
\Zp\Zp\Zp\Zm\Zz\\
\Zp\Za\ZA\Zz\Zm\\
\Zp\ZA\Zz\Zb\ZB\\
\Zp\Zz\Za\ZB\Zb\\
\Zz\Zp\ZA\Zb\Zb\\
\end{array}\right)$ & $
\left(\begin{array}{c}
\Zp\Zp\Zm\Zp\Zz\\
\Zp\Za\ZB\Zz\Zm\\
\Zp\ZA\Zz\Za\ZB\\
\Zp\Zz\Zb\ZA\Zb\\
\Zz\Zp\ZB\Za\Zb\\
\end{array}\right)$ & $
\left(\begin{array}{c}
\Zp\ZB\Zz\Za\ZB\\
\Zp\Zm\Zp\Zp\Zz\\
\Zp\Zb\ZA\Zz\Zm\\
\Zp\Zz\Za\ZA\Zb\\
\Zz\Zp\ZB\Zb\Za\\
\end{array}\right)$ & $
\left(\begin{array}{c}
\Zp\ZB\Zz\Zb\ZA\\
\Zp\Zm\Zm\Zm\Zz\\
\Zp\Zb\ZB\Zz\Zp\\
\Zp\Zz\Zb\ZB\Za\\
\Zz\Zp\ZA\Za\Zb\\
\end{array}\right)
$\\
\bottomrule
\end{tabular}
\end{table}

\begin{table}[H]\centering\caption{20 mutually unbiased weighing matrices of order 6 and weight 4, $UW(6,4)$.}
\begin{tabular}{@{}cccc@{}}\toprule\label{table:UW6_4}
$
\left(\begin{array}{c}
\Zp\Zp\Zp\Zp\Zz\Zz\\
\Zp\Zp\Zm\Zm\Zz\Zz\\
\Zp\Zm\Zz\Zz\Zp\Zp\\
\Zp\Zm\Zz\Zz\Zm\Zm\\
\Zz\Zz\Zp\Zm\Zp\Zm\\
\Zz\Zz\Zp\Zm\Zm\Zp
\end{array}\right)$ & $
\left(\begin{array}{c}
\Zp\Zp\Zp\Zm\Zz\Zz\\
\Zp\Zp\Zm\Zp\Zz\Zz\\
\Zp\Zm\Zz\Zz\Zp\Zm\\
\Zp\Zm\Zz\Zz\Zm\Zp\\
\Zz\Zz\Zp\Zp\Zp\Zp\\
\Zz\Zz\Zp\Zp\Zm\Zm
\end{array}\right)$ & $
\left(\begin{array}{c}
\Zp\Zp\Zz\Zz\Zp\Zp\\
\Zp\Zp\Zz\Zz\Zm\Zm\\
\Zp\Zm\Zp\Zm\Zz\Zz\\
\Zp\Zm\Zm\Zp\Zz\Zz\\
\Zz\Zz\Zp\Zp\Zp\Zm\\
\Zz\Zz\Zp\Zp\Zm\Zp
\end{array}\right)$ & $
\left(\begin{array}{c}
\Zp\Zp\Zz\Zz\Zp\Zm\\
\Zp\Zp\Zz\Zz\Zm\Zp\\
\Zp\Zm\Zp\Zp\Zz\Zz\\
\Zp\Zm\Zm\Zm\Zz\Zz\\
\Zz\Zz\Zp\Zm\Zp\Zp\\
\Zz\Zz\Zp\Zm\Zm\Zm
\end{array}\right)$ \\$
\left(\begin{array}{c}
\Zp\ZB\Za\Zz\Za\Zz\\
\Zp\ZA\Za\Zz\Zb\Zz\\
\Zp\Zz\Zb\ZB\Zz\Zp\\
\Zp\Zz\Zb\ZA\Zz\Zm\\
\Zz\Zp\Zz\Zp\ZA\Za\\
\Zz\Zp\Zz\Zm\ZA\Zb
\end{array}\right)$ & $
\left(\begin{array}{c}
\Zp\ZB\Za\Zz\Zb\Zz\\
\Zp\ZA\Za\Zz\Za\Zz\\
\Zp\Zz\Zb\ZB\Zz\Zm\\
\Zp\Zz\Zb\ZA\Zz\Zp\\
\Zz\Zp\Zz\Zp\ZB\Zb\\
\Zz\Zp\Zz\Zm\ZB\Za
\end{array}\right)$ & $
\left(\begin{array}{c}
\Zp\ZB\Zb\Zz\Za\Zz\\
\Zp\ZA\Zb\Zz\Zb\Zz\\
\Zp\Zz\Za\ZB\Zz\Zm\\
\Zp\Zz\Za\ZA\Zz\Zp\\
\Zz\Zp\Zz\Zp\ZA\Zb\\
\Zz\Zp\Zz\Zm\ZA\Za
\end{array}\right)$ & $
\left(\begin{array}{c}
\Zp\ZB\Zb\Zz\Zb\Zz\\
\Zp\ZA\Zb\Zz\Za\Zz\\
\Zp\Zz\Za\ZB\Zz\Zp\\
\Zp\Zz\Za\ZA\Zz\Zm\\
\Zz\Zp\Zz\Zp\ZB\Za\\
\Zz\Zp\Zz\Zm\ZB\Zb
\end{array}\right)$ \\$
\left(\begin{array}{c}
\Zp\ZB\Zz\Za\Zz\Za\\
\Zp\ZA\Zz\Za\Zz\Zb\\
\Zp\Zz\ZB\Zb\Zp\Zz\\
\Zp\Zz\ZA\Zb\Zm\Zz\\
\Zz\Zp\Zp\Zz\Za\ZA\\
\Zz\Zp\Zm\Zz\Zb\ZA
\end{array}\right)$ & $
\left(\begin{array}{c}
\Zp\ZB\Zz\Za\Zz\Zb\\
\Zp\ZA\Zz\Za\Zz\Za\\
\Zp\Zz\ZB\Zb\Zm\Zz\\
\Zp\Zz\ZA\Zb\Zp\Zz\\
\Zz\Zp\Zp\Zz\Zb\ZB\\
\Zz\Zp\Zm\Zz\Za\ZB
\end{array}\right)$ & $
\left(\begin{array}{c}
\Zp\ZB\Zz\Zb\Zz\Za\\
\Zp\ZA\Zz\Zb\Zz\Zb\\
\Zp\Zz\ZB\Za\Zm\Zz\\
\Zp\Zz\ZA\Za\Zp\Zz\\
\Zz\Zp\Zp\Zz\Zb\ZA\\
\Zz\Zp\Zm\Zz\Za\ZA
\end{array}\right)$ & $
\left(\begin{array}{c}
\Zp\ZB\Zz\Zb\Zz\Zb\\
\Zp\ZA\Zz\Zb\Zz\Za\\
\Zp\Zz\ZB\Za\Zp\Zz\\
\Zp\Zz\ZA\Za\Zm\Zz\\
\Zz\Zp\Zp\Zz\Za\ZB\\
\Zz\Zp\Zm\Zz\Zb\ZB
\end{array}\right)$ \\$
\left(\begin{array}{c}
\Zp\Za\ZB\Zz\Zz\ZB\\
\Zp\Zb\ZA\Zz\Zz\ZB\\
\Zp\Zz\Zz\Zp\Za\ZA\\
\Zp\Zz\Zz\Zm\Zb\ZA\\
\Zz\Zp\Za\ZB\Zp\Zz\\
\Zz\Zp\Za\ZA\Zm\Zz
\end{array}\right)$ & $
\left(\begin{array}{c}
\Zp\Za\ZB\Zz\Zz\ZA\\
\Zp\Zb\ZA\Zz\Zz\ZA\\
\Zp\Zz\Zz\Zp\Zb\ZB\\
\Zp\Zz\Zz\Zm\Za\ZB\\
\Zz\Zp\Za\ZB\Zm\Zz\\
\Zz\Zp\Za\ZA\Zp\Zz
\end{array}\right)$ & $
\left(\begin{array}{c}
\Zp\Za\ZA\Zz\Zz\ZB\\
\Zp\Zb\ZB\Zz\Zz\ZB\\
\Zp\Zz\Zz\Zp\Zb\ZA\\
\Zp\Zz\Zz\Zm\Za\ZA\\
\Zz\Zp\Zb\ZB\Zm\Zz\\
\Zz\Zp\Zb\ZA\Zp\Zz
\end{array}\right)$ & $
\left(\begin{array}{c}
\Zp\Za\ZA\Zz\Zz\ZA\\
\Zp\Zb\ZB\Zz\Zz\ZA\\
\Zp\Zz\Zz\Zp\Za\ZB\\
\Zp\Zz\Zz\Zm\Zb\ZB\\
\Zz\Zp\Zb\ZB\Zp\Zz\\
\Zz\Zp\Zb\ZA\Zm\Zz
\end{array}\right)$ \\$
\left(\begin{array}{c}
\Zp\Za\Zz\ZB\ZB\Zz\\
\Zp\Zb\Zz\ZA\ZB\Zz\\
\Zp\Zz\Zp\Zz\ZA\Za\\
\Zp\Zz\Zm\Zz\ZA\Zb\\
\Zz\Zp\ZB\Za\Zz\Zp\\
\Zz\Zp\ZA\Za\Zz\Zm
\end{array}\right)$ & $
\left(\begin{array}{c}
\Zp\Za\Zz\ZB\ZA\Zz\\
\Zp\Zb\Zz\ZA\ZA\Zz\\
\Zp\Zz\Zp\Zz\ZB\Zb\\
\Zp\Zz\Zm\Zz\ZB\Za\\
\Zz\Zp\ZB\Za\Zz\Zm\\
\Zz\Zp\ZA\Za\Zz\Zp
\end{array}\right)$ & $
\left(\begin{array}{c}
\Zp\Za\Zz\ZA\ZB\Zz\\
\Zp\Zb\Zz\ZB\ZB\Zz\\
\Zp\Zz\Zp\Zz\ZA\Zb\\
\Zp\Zz\Zm\Zz\ZA\Za\\
\Zz\Zp\ZB\Zb\Zz\Zm\\
\Zz\Zp\ZA\Zb\Zz\Zp
\end{array}\right)$ & $
\left(\begin{array}{c}
\Zp\Za\Zz\ZA\ZA\Zz\\
\Zp\Zb\Zz\ZB\ZA\Zz\\
\Zp\Zz\Zp\Zz\ZB\Za\\
\Zp\Zz\Zm\Zz\ZB\Zb\\
\Zz\Zp\ZB\Zb\Zz\Zp\\
\Zz\Zp\ZA\Zb\Zz\Zm
\end{array}\right)
$\\
\bottomrule
\end{tabular}
\end{table}

\begin{table}[H]\centering\caption{8 mutually unbiased real weighing matrices of order 7 and weight 4, $W(7,4)$.}
\begin{tabular}{@{}cccc@{}}\toprule\label{table:W7_4}
$\left(
\begin{array}{c}
\Zp\Zp\Zp\Zp\Zz\Zz\Zz\\
\Zp\Zm\Zz\Zz\Zp\Zp\Zz\\
\Zp\Zz\Zm\Zz\Zm\Zz\Zp\\
\Zp\Zz\Zz\Zm\Zz\Zm\Zm\\
\Zz\Zp\Zm\Zz\Zz\Zp\Zm\\
\Zz\Zp\Zz\Zm\Zp\Zz\Zp\\
\Zz\Zz\Zp\Zm\Zm\Zp\Zz\\
\end{array}
\right)$ &
$\left(
\begin{array}{c}
\Zp\Zp\Zm\Zm\Zz\Zz\Zz\\
\Zp\Zm\Zz\Zz\Zm\Zp\Zz\\
\Zp\Zz\Zp\Zz\Zp\Zz\Zp\\
\Zp\Zz\Zz\Zp\Zz\Zm\Zm\\
\Zz\Zp\Zp\Zz\Zz\Zp\Zm\\
\Zz\Zp\Zz\Zp\Zm\Zz\Zp\\
\Zz\Zz\Zp\Zm\Zm\Zm\Zz\\
\end{array}
\right)$ &
$\left(
\begin{array}{c}
\Zp\Zp\Zz\Zz\Zm\Zm\Zz\\
\Zp\Zm\Zp\Zm\Zz\Zz\Zz\\
\Zp\Zz\Zm\Zz\Zp\Zz\Zp\\
\Zp\Zz\Zz\Zp\Zz\Zp\Zm\\
\Zz\Zp\Zp\Zz\Zz\Zp\Zp\\
\Zz\Zp\Zz\Zm\Zp\Zz\Zm\\
\Zz\Zz\Zp\Zp\Zp\Zm\Zz\\
\end{array}
\right)$ &
$\left(
\begin{array}{c}
\Zp\Zp\Zz\Zz\Zp\Zp\Zz\\
\Zp\Zm\Zm\Zm\Zz\Zz\Zz\\
\Zp\Zz\Zp\Zz\Zm\Zz\Zm\\
\Zp\Zz\Zz\Zp\Zz\Zm\Zp\\
\Zz\Zp\Zm\Zz\Zz\Zm\Zm\\
\Zz\Zp\Zz\Zm\Zm\Zz\Zp\\
\Zz\Zz\Zp\Zm\Zp\Zm\Zz\\
\end{array}
\right)$
\\
$\left(
\begin{array}{c}
\Zp\Zp\Zp\Zm\Zz\Zz\Zz\\
\Zp\Zm\Zz\Zz\Zp\Zm\Zz\\
\Zp\Zz\Zm\Zz\Zm\Zz\Zm\\
\Zp\Zz\Zz\Zp\Zz\Zp\Zp\\
\Zz\Zp\Zm\Zz\Zz\Zm\Zp\\
\Zz\Zp\Zz\Zp\Zp\Zz\Zm\\
\Zz\Zz\Zp\Zp\Zm\Zm\Zz\\
\end{array}
\right)$ &
$\left(
\begin{array}{c}
\Zp\Zp\Zz\Zz\Zp\Zm\Zz\\
\Zp\Zm\Zm\Zp\Zz\Zz\Zz\\
\Zp\Zz\Zp\Zz\Zm\Zz\Zp\\
\Zp\Zz\Zz\Zm\Zz\Zp\Zm\\
\Zz\Zp\Zm\Zz\Zz\Zp\Zp\\
\Zz\Zp\Zz\Zp\Zm\Zz\Zm\\
\Zz\Zz\Zp\Zp\Zp\Zp\Zz\\
\end{array}
\right)$ &
$\left(
\begin{array}{c}
\Zp\Zp\Zm\Zp\Zz\Zz\Zz\\
\Zp\Zm\Zz\Zz\Zm\Zm\Zz\\
\Zp\Zz\Zp\Zz\Zp\Zz\Zm\\
\Zp\Zz\Zz\Zm\Zz\Zp\Zp\\
\Zz\Zp\Zp\Zz\Zz\Zm\Zp\\
\Zz\Zp\Zz\Zm\Zm\Zz\Zm\\
\Zz\Zz\Zp\Zp\Zm\Zp\Zz\\
\end{array}
\right)$ &
$\left(
\begin{array}{c}
\Zp\Zp\Zz\Zz\Zm\Zp\Zz\\
\Zp\Zm\Zp\Zp\Zz\Zz\Zz\\
\Zp\Zz\Zm\Zz\Zp\Zz\Zm\\
\Zp\Zz\Zz\Zm\Zz\Zm\Zp\\
\Zz\Zp\Zp\Zz\Zz\Zm\Zm\\
\Zz\Zp\Zz\Zp\Zp\Zz\Zp\\
\Zz\Zz\Zp\Zm\Zp\Zp\Zz
\end{array}
\right)$\\
\bottomrule
\end{tabular}
\end{table}

\begin{table}[H]\centering\caption{14 mutually unbiased real weighing matrices of order 8 and weight 4, $W(8,4)$.}
\begin{tabular}{@{}cccc@{}}\toprule\label{table:W8_4}
$\left(
\begin{array}{c}
\Zp\Zp\Zp\Zp\Zz\Zz\Zz\Zz\\
\Zp\Zp\Zm\Zm\Zz\Zz\Zz\Zz\\
\Zp\Zm\Zp\Zm\Zz\Zz\Zz\Zz\\
\Zp\Zm\Zm\Zp\Zz\Zz\Zz\Zz\\
\Zz\Zz\Zz\Zz\Zp\Zp\Zp\Zp\\
\Zz\Zz\Zz\Zz\Zp\Zp\Zm\Zm\\
\Zz\Zz\Zz\Zz\Zp\Zm\Zp\Zm\\
\Zz\Zz\Zz\Zz\Zp\Zm\Zm\Zp\\
\end{array}
\right)$ &
$\left(
\begin{array}{c}
\Zp\Zp\Zp\Zm\Zz\Zz\Zz\Zz\\
\Zp\Zp\Zm\Zp\Zz\Zz\Zz\Zz\\
\Zp\Zm\Zz\Zz\Zp\Zm\Zz\Zz\\
\Zp\Zm\Zz\Zz\Zm\Zp\Zz\Zz\\
\Zz\Zz\Zp\Zp\Zz\Zz\Zp\Zp\\
\Zz\Zz\Zp\Zp\Zz\Zz\Zm\Zm\\
\Zz\Zz\Zz\Zz\Zp\Zp\Zp\Zm\\
\Zz\Zz\Zz\Zz\Zp\Zp\Zm\Zp\\
\end{array}
\right)$ &
$\left(
\begin{array}{c}
\Zp\Zp\Zz\Zz\Zp\Zp\Zz\Zz\\
\Zp\Zp\Zz\Zz\Zm\Zm\Zz\Zz\\
\Zp\Zm\Zp\Zp\Zz\Zz\Zz\Zz\\
\Zp\Zm\Zm\Zm\Zz\Zz\Zz\Zz\\
\Zz\Zz\Zp\Zm\Zz\Zz\Zp\Zm\\
\Zz\Zz\Zp\Zm\Zz\Zz\Zm\Zp\\
\Zz\Zz\Zz\Zz\Zp\Zm\Zp\Zp\\
\Zz\Zz\Zz\Zz\Zp\Zm\Zm\Zm\\
\end{array}
\right)$ &
$\left(
\begin{array}{c}
\Zp\Zz\Zz\Zp\Zz\Zp\Zm\Zz\\
\Zp\Zz\Zz\Zp\Zz\Zm\Zp\Zz\\
\Zp\Zz\Zz\Zm\Zp\Zz\Zz\Zm\\
\Zp\Zz\Zz\Zm\Zm\Zz\Zz\Zp\\
\Zz\Zp\Zp\Zz\Zp\Zz\Zz\Zp\\
\Zz\Zp\Zp\Zz\Zm\Zz\Zz\Zm\\
\Zz\Zp\Zm\Zz\Zz\Zp\Zp\Zz\\
\Zz\Zp\Zm\Zz\Zz\Zm\Zm\Zz\\
\end{array}
\right)$
\\
$\left(
\begin{array}{c}
\Zp\Zp\Zz\Zz\Zp\Zm\Zz\Zz\\
\Zp\Zp\Zz\Zz\Zm\Zp\Zz\Zz\\
\Zp\Zm\Zz\Zz\Zz\Zz\Zp\Zp\\
\Zp\Zm\Zz\Zz\Zz\Zz\Zm\Zm\\
\Zz\Zz\Zp\Zp\Zz\Zz\Zp\Zm\\
\Zz\Zz\Zp\Zp\Zz\Zz\Zm\Zp\\
\Zz\Zz\Zp\Zm\Zp\Zp\Zz\Zz\\
\Zz\Zz\Zp\Zm\Zm\Zm\Zz\Zz\\
\end{array}
\right)$ &
$\left(
\begin{array}{c}
\Zp\Zz\Zz\Zp\Zp\Zz\Zz\Zp\\
\Zp\Zz\Zz\Zp\Zm\Zz\Zz\Zm\\
\Zp\Zz\Zz\Zm\Zz\Zp\Zp\Zz\\
\Zp\Zz\Zz\Zm\Zz\Zm\Zm\Zz\\
\Zz\Zp\Zp\Zz\Zz\Zp\Zm\Zz\\
\Zz\Zp\Zp\Zz\Zz\Zm\Zp\Zz\\
\Zz\Zp\Zm\Zz\Zp\Zz\Zz\Zm\\
\Zz\Zp\Zm\Zz\Zm\Zz\Zz\Zp\\
\end{array}
\right)$ &
$\left(
\begin{array}{c}
\Zp\Zp\Zz\Zz\Zz\Zz\Zp\Zm\\
\Zp\Zp\Zz\Zz\Zz\Zz\Zm\Zp\\
\Zp\Zm\Zz\Zz\Zp\Zp\Zz\Zz\\
\Zp\Zm\Zz\Zz\Zm\Zm\Zz\Zz\\
\Zz\Zz\Zp\Zp\Zp\Zm\Zz\Zz\\
\Zz\Zz\Zp\Zp\Zm\Zp\Zz\Zz\\
\Zz\Zz\Zp\Zm\Zz\Zz\Zp\Zp\\
\Zz\Zz\Zp\Zm\Zz\Zz\Zm\Zm\\
\end{array}
\right)$ &
$\left(
\begin{array}{c}
\Zp\Zz\Zz\Zp\Zp\Zz\Zz\Zm\\
\Zp\Zz\Zz\Zp\Zm\Zz\Zz\Zp\\
\Zp\Zz\Zz\Zm\Zp\Zz\Zz\Zp\\
\Zp\Zz\Zz\Zm\Zm\Zz\Zz\Zm\\
\Zz\Zp\Zp\Zz\Zz\Zp\Zp\Zz\\
\Zz\Zp\Zp\Zz\Zz\Zm\Zm\Zz\\
\Zz\Zp\Zm\Zz\Zz\Zp\Zm\Zz\\
\Zz\Zp\Zm\Zz\Zz\Zm\Zp\Zz\\
\end{array}
\right)$
\\
$\left(
\begin{array}{c}
\Zp\Zz\Zz\Zp\Zz\Zp\Zp\Zz\\
\Zp\Zz\Zz\Zp\Zz\Zm\Zm\Zz\\
\Zp\Zz\Zz\Zm\Zz\Zp\Zm\Zz\\
\Zp\Zz\Zz\Zm\Zz\Zm\Zp\Zz\\
\Zz\Zp\Zp\Zz\Zp\Zz\Zz\Zm\\
\Zz\Zp\Zp\Zz\Zm\Zz\Zz\Zp\\
\Zz\Zp\Zm\Zz\Zp\Zz\Zz\Zp\\
\Zz\Zp\Zm\Zz\Zm\Zz\Zz\Zm\\
\end{array}
\right)$ &
$\left(
\begin{array}{c}
\Zp\Zz\Zp\Zz\Zz\Zp\Zz\Zm\\
\Zp\Zz\Zp\Zz\Zz\Zm\Zz\Zp\\
\Zp\Zz\Zm\Zz\Zp\Zz\Zp\Zz\\
\Zp\Zz\Zm\Zz\Zm\Zz\Zm\Zz\\
\Zz\Zp\Zz\Zp\Zp\Zz\Zm\Zz\\
\Zz\Zp\Zz\Zp\Zm\Zz\Zp\Zz\\
\Zz\Zp\Zz\Zm\Zz\Zp\Zz\Zp\\
\Zz\Zp\Zz\Zm\Zz\Zm\Zz\Zm\\
\end{array}
\right)$ &
$\left(
\begin{array}{c}
\Zp\Zz\Zp\Zz\Zp\Zz\Zp\Zz\\
\Zp\Zz\Zp\Zz\Zm\Zz\Zm\Zz\\
\Zp\Zz\Zm\Zz\Zp\Zz\Zm\Zz\\
\Zp\Zz\Zm\Zz\Zm\Zz\Zp\Zz\\
\Zz\Zp\Zz\Zp\Zz\Zp\Zz\Zp\\
\Zz\Zp\Zz\Zp\Zz\Zm\Zz\Zm\\
\Zz\Zp\Zz\Zm\Zz\Zp\Zz\Zm\\
\Zz\Zp\Zz\Zm\Zz\Zm\Zz\Zp\\
\end{array}
\right)$ &
$\left(
\begin{array}{c}
\Zp\Zp\Zz\Zz\Zz\Zz\Zp\Zp\\
\Zp\Zm\Zz\Zz\Zz\Zz\Zm\Zp\\
\Zp\Zz\Zp\Zz\Zz\Zm\Zz\Zm\\
\Zp\Zz\Zm\Zz\Zz\Zp\Zz\Zm\\
\Zz\Zp\Zz\Zp\Zm\Zz\Zm\Zz\\
\Zz\Zp\Zz\Zm\Zp\Zz\Zm\Zz\\
\Zz\Zz\Zp\Zp\Zp\Zp\Zz\Zz\\
\Zz\Zz\Zp\Zm\Zm\Zp\Zz\Zz\\
\end{array}
\right)$
\\
$\left(
\begin{array}{c}
\Zp\Zp\Zz\Zz\Zz\Zz\Zm\Zm\\
\Zp\Zm\Zz\Zz\Zz\Zz\Zp\Zm\\
\Zp\Zz\Zp\Zz\Zz\Zp\Zz\Zp\\
\Zp\Zz\Zm\Zz\Zz\Zm\Zz\Zp\\
\Zz\Zp\Zz\Zp\Zp\Zz\Zp\Zz\\
\Zz\Zp\Zz\Zm\Zm\Zz\Zp\Zz\\
\Zz\Zz\Zp\Zp\Zm\Zm\Zz\Zz\\
\Zz\Zz\Zp\Zm\Zp\Zm\Zz\Zz\\
\end{array}
\right)$ &
$\left(
\begin{array}{c}
\Zp\Zz\Zp\Zz\Zp\Zz\Zm\Zz\\
\Zp\Zz\Zp\Zz\Zm\Zz\Zp\Zz\\
\Zp\Zz\Zm\Zz\Zz\Zp\Zz\Zp\\
\Zp\Zz\Zm\Zz\Zz\Zm\Zz\Zm\\
\Zz\Zp\Zz\Zp\Zz\Zp\Zz\Zm\\
\Zz\Zp\Zz\Zp\Zz\Zm\Zz\Zp\\
\Zz\Zp\Zz\Zm\Zp\Zz\Zp\Zz\\
\Zz\Zp\Zz\Zm\Zm\Zz\Zm\Zz
\end{array}
\right)$\\
\bottomrule
\end{tabular}
\end{table}

\section{Hadamard matrices of order 32}\label{app:H32}
In Tables \ref{table:H32_1},\ref{table:H32_2},\ref{table:H32_3} and \ref{table:H32_4}, we show the partition of the $32^2$ vectors into $32$ Hadamard matrices of order $32$ (denoted $H_1,H_2,\cdots,H_{32}$). Each section represents one Hadamard matrix, and each hexadecimal number represents one row of the matrix (where each digit represents four entries).

\begin{table}\centering\caption{$H_{1}$ through $H_{8}$}
\begin{tabular}{@{}cccccccc@{}}\toprule\label{table:H32_1}
\texttt{00000000} & \texttt{4259F1BA} & \texttt{203AEEB5} & \texttt{50967C6E} & \texttt{59F1BA84} & \texttt{47FC04A7} & \texttt{4E9BC24D} & \texttt{4B3E3750} \\
\texttt{62631F0F} & \texttt{7C6EA12C} & \texttt{1E0DBE23} & \texttt{259F1BA8} & \texttt{32F56361} & \texttt{750967C6} & \texttt{176A78C9} & \texttt{67C6EA12} \\
\texttt{6EA12CF8} & \texttt{70AC92DB} & \texttt{55338973} & \texttt{0CC233F7} & \texttt{12CF8DD4} & \texttt{05A5F51D} & \texttt{3B92A58B} & \texttt{79CB5431} \\
\texttt{1BA84B3E} & \texttt{5C544F99} & \texttt{6B04D9E5} & \texttt{3E375096} & \texttt{2CF8DD42} & \texttt{0967C6EA} & \texttt{3750967C} & \texttt{295D285F} \\

\midrule
\texttt{6EF49ECD} & \texttt{2D2FA8E1} & \texttt{755CD5F3} & \texttt{1BFDF90B} & \texttt{31A55E78} & \texttt{7C3B1319} & \texttt{3FE02535} & \texttt{5826CF27} \\
\texttt{727E6854} & \texttt{4DCBFF54} & \texttt{5663B46A} & \texttt{7B19AEBE} & \texttt{129A3FE1} & \texttt{0055B235} & \texttt{24486E0B} & \texttt{44AC39BE} \\
\texttt{0E10C978} & \texttt{38C29892} & \texttt{4AE942F3} & \texttt{15B88246} & \texttt{438E8419} & \texttt{514109CD} & \texttt{67935827} & \texttt{2A0D1546} \\
\texttt{69D6236A} & \texttt{3687E3DF} & \texttt{093274DF} & \texttt{1CDF44AC} & \texttt{60B1E580} & \texttt{5F047280} & \texttt{236AD3AC} & \texttt{07770F92} \\

\midrule
\texttt{477DB9D5} & \texttt{1F0EC4C7} & \texttt{6A07A301} & \texttt{182C7960} & \texttt{0384325E} & \texttt{5CD5F2EB} & \texttt{5BF74F4C} & \texttt{529089A6} \\
\texttt{636065EB} & \texttt{114BBF8A} & \texttt{7FEA9372} & \texttt{2EFE288A} & \texttt{405F0472} & \texttt{71AFE83F} & \texttt{4E1A7F3F} & \texttt{2799EE60} \\
\texttt{0AE3F4B4} & \texttt{0DC14913} & \texttt{355663B4} & \texttt{20BB53C7} & \texttt{78C82ED5} & \texttt{29DC952D} & \texttt{768D5598} & \texttt{3274DE13} \\
\texttt{1669022D} & \texttt{3C31A55E} & \texttt{4938C298} & \texttt{04A68FF9} & \texttt{6D251EA6} & \texttt{6442D84C} & \texttt{55B23401} & \texttt{3B1318F9} \\

\midrule
\texttt{050E9174} & \texttt{10E3A107} & \texttt{19D1D5D8} & \texttt{475760CE} & \texttt{7935826D} & \texttt{20C438E9} & \texttt{6BAFBD8C} & \texttt{629DC953} \\
\texttt{2E8143A4} & \texttt{3529089A} & \texttt{4E651411} & \texttt{7770F920} & \texttt{52BA50BD} & \texttt{02799EE6} & \texttt{3C1B7C45} & \texttt{40206F5C} \\
\texttt{5CFF2BF0} & \texttt{7E428DFF} & \texttt{27B3377B} & \texttt{29F64C36} & \texttt{65EAC6C1} & \texttt{1EA6DA4A} & \texttt{0B4BEA39} & \texttt{325E0708} \\
\texttt{3B6C73D7} & \texttt{5B882462} & \texttt{7007F6B2} & \texttt{49121B83} & \texttt{6CD8B21E} & \texttt{1794AE95} & \texttt{0C3CE5AB} & \texttt{55CD5F2F} \\

\midrule
\texttt{7357CBAB} & \texttt{66BAFBD8} & \texttt{7475760C} & \texttt{017C11CA} & \texttt{5FD07DC7} & \texttt{2AA6712F} & \texttt{3F342A72} & \texttt{5FAF16E9} \\
\texttt{14EE4A97} & \texttt{3F4B415C} & \texttt{61E72D51} & \texttt{4D1FF013} & \texttt{0621C743} & \texttt{381697D5} & \texttt{4A3D4DB4} & \texttt{149121B9} \\
\texttt{7328A085} & \texttt{740A1D22} & \texttt{01037AE4} & \texttt{66C590F6} & \texttt{2D84CC88} & \texttt{58F2C060} & \texttt{065EAC6D} & \texttt{4A42269A} \\
\texttt{13CCF730} & \texttt{2DFBA7A6} & \texttt{3869FCFB} & \texttt{4D609B3D} & \texttt{2AD91A01} & \texttt{588DAB4E} & \texttt{13B39C1E} & \texttt{6198467F} \\

\midrule
\texttt{76595ADF} & \texttt{03503D19} & \texttt{53C64177} & \texttt{1CF59DB7} & \texttt{0D154654} & \texttt{0E3A1063} & \texttt{4C63E1D9} & \texttt{1FDACB80} \\
\texttt{288A5DFC} & \texttt{5EAC6C0D} & \texttt{372FFD52} & \texttt{4109CCA3} & \texttt{12B0E6FA} & \texttt{6496D70B} & \texttt{119FB0CD} & \texttt{4F4CB7EE} \\
\texttt{396A861F} & \texttt{007F6B2E} & \texttt{50E91740} & \texttt{69FCFA71} & \texttt{781C2192} & \texttt{5D833A3A} & \texttt{3400AB65} & \texttt{42269A94} \\
\texttt{25E07086} & \texttt{75760CE8} & \texttt{3A45D028} & \texttt{2BA50BCB} & \texttt{26CF26B1} & \texttt{7B3377A5} & \texttt{67B9813C} & \texttt{6AD3AC46} \\

\midrule
\texttt{7A4F666F} & \texttt{0F4601A9} & \texttt{5195068A} & \texttt{7A300D41} & \texttt{56B7BB2D} & \texttt{23BEDCEB} & \texttt{44075DD7} & \texttt{3171513F} \\
\texttt{24E30A62} & \texttt{56C8D003} & \texttt{0864BC0E} & \texttt{435A8B5E} & \texttt{1DF6E753} & \texttt{1AAB31DA} & \texttt{23C1B7C5} & \texttt{6880EBBB} \\
\texttt{0F396A87} & \texttt{7D6DDBC8} & \texttt{7D12B0E6} & \texttt{447836F9} & \texttt{4325E070} & \texttt{310E3A11} & \texttt{68FF8095} & \texttt{362C87B6} \\
\texttt{1AD45AF4} & \texttt{3653EC98} & \texttt{6FDD3D32} & \texttt{249C614C} & \texttt{1D898C7D} & \texttt{51EA6DA4} & \texttt{081BD720} & \texttt{6FA2561C} \\

\midrule
\texttt{0BCA574B} & \texttt{621C7421} & \texttt{70864BC0} & \texttt{4993A6F1} & \texttt{27328A09} & \texttt{2E00FED6} & \texttt{523BEDCF} & \texttt{32DFBA7A} \\
\texttt{5C7E9682} & \texttt{7EC3308D} & \texttt{77F14452} & \texttt{3C9AC137} & \texttt{171513E7} & \texttt{10621C75} & \texttt{4EE4A963} & \texttt{1E276738} \\
\texttt{058F2C06} & \texttt{195068AA} & \texttt{02F82394} & \texttt{3BEDCEA5} & \texttt{0CBD58D9} & \texttt{2045859B} & \texttt{5B099910} & \texttt{79B43F1F} \\
\texttt{47D6DDBC} & \texttt{656B7BB3} & \texttt{554CE25D} & \texttt{6C590F6C} & \texttt{35A8B5E8} & \texttt{6B2E00FE} & \texttt{40A1D22E} & \texttt{2977F144} \\

\bottomrule
\end{tabular}
\end{table}

\begin{table}\centering\caption{$H_{9}$ through $H_{16}$}
\begin{tabular}{@{}cccccccc@{}}\toprule\label{table:H32_2}
\texttt{09B3C9AD} & \texttt{7B9813CC} & \texttt{7DC6BFA1} & \texttt{2C2CD205} & \texttt{2A727E68} & \texttt{513E62E3} & \texttt{74A1794B} & \texttt{13679359} \\
\texttt{00D40F47} & \texttt{254B14EF} & \texttt{2315B882} & \texttt{3F9F4E1B} & \texttt{4DB4947A} & \texttt{44D35290} & \texttt{0FED65C0} & \texttt{1A0055B3} \\
\texttt{614C4938} & \texttt{5859A409} & \texttt{72FFD526} & \texttt{6712E555} & \texttt{068AA32A} & \texttt{6E7523BF} & \texttt{36F888F1} & \texttt{1C5EF9DE} \\
\texttt{30A6249C} & \texttt{5760CE8E} & \texttt{682B8FD2} & \texttt{428DFEFD} & \texttt{5E070864} & \texttt{4BEA3817} & \texttt{39C1E276} & \texttt{15393F34} \\

\midrule
\texttt{1853124E} & \texttt{16166903} & \texttt{146FF7E5} & \texttt{4680156D} & \texttt{78B745FB} & \texttt{0FC7BCDB} & \texttt{241DDC3E} & \texttt{1A2A8CA8} \\
\texttt{266442D8} & \texttt{2A58A773} & \texttt{6F23EB6E} & \texttt{631F0EC5} & \texttt{3DCC09E6} & \texttt{3FB59700} & \texttt{48C56E20} & \texttt{4ABCF0C6} \\
\texttt{61669023} & \texttt{6D5A7588} & \texttt{0182C796} & \texttt{28213995} & \texttt{0DBE223D} & \texttt{536D251E} & \texttt{5F51C0B5} & \texttt{5114BBF8} \\
\texttt{76F23EB6} & \texttt{03FB5970} & \texttt{5D285E53} & \texttt{338972AB} & \texttt{31F0EC4D} & \texttt{7ACEDB1D} & \texttt{748BA050} & \texttt{44F98B8B} \\

\midrule
\texttt{1FA5A0AE} & \texttt{2FD78B75} & \texttt{111E0DBF} & \texttt{347FC04B} & \texttt{5347FC05} & \texttt{789D9CE0} & \texttt{46541A2A} & \texttt{040DEB90} \\
\texttt{48116167} & \texttt{216C2664} & \texttt{6D70AC93} & \texttt{5D028748} & \texttt{3AC46D5A} & \texttt{1F5B76F2} & \texttt{6D8E7ACF} & \texttt{46AACC76} \\
\texttt{0AB64681} & \texttt{0A4890DD} & \texttt{53B92A59} & \texttt{48EFB73B} & \texttt{2F295D29} & \texttt{76D8E7AD} & \texttt{3A3ABB06} & \texttt{6335D7DE} \\
\texttt{2192F038} & \texttt{762631F1} & \texttt{04F33DCC} & \texttt{78634ABC} & \texttt{34811617} & \texttt{63CB0182} & \texttt{5DFC5114} & \texttt{11E0DBE3} \\

\midrule
\texttt{6C0CBD59} & \texttt{3BB87C90} & \texttt{1740A1D2} & \texttt{2E554CE3} & \texttt{47836F89} & \texttt{05DA9E33} & \texttt{35FD07DD} & \texttt{40F4601B} \\
\texttt{526E5FFA} & \texttt{201037AE} & \texttt{77A4F667} & \texttt{1E72D50D} & \texttt{2767383C} & \texttt{29224371} & \texttt{49C614C4} & \texttt{6249C614} \\
\texttt{4EB11B56} & \texttt{3CCF7302} & \texttt{5C2B24B7} & \texttt{1037AE40} & \texttt{7E9682B8} & \texttt{5B5C2B25} & \texttt{6B7BB2CB} & \texttt{02AD91A1} \\
\texttt{0B9FE57E} & \texttt{0CE8EAEC} & \texttt{653EC986} & \texttt{79E18D2A} & \texttt{70D3F9F5} & \texttt{1905DA9F} & \texttt{328A084F} & \texttt{55195068} \\

\midrule
\texttt{5017C11C} & \texttt{3D32DFBA} & \texttt{11CA02F8} & \texttt{7CEF1C5E} & \texttt{415C7E96} & \texttt{784993A7} & \texttt{0427328B} & \texttt{37D12B0E} \\
\texttt{69022C2D} & \texttt{4F1905DB} & \texttt{5E52BA51} & \texttt{63E1D899} & \texttt{760CE8EA} & \texttt{5AF435A8} & \texttt{0081BD72} & \texttt{39945043} \\
\texttt{28DFEFC9} & \texttt{223C1B7D} & \texttt{1B29F64C} & \texttt{269A9484} & \texttt{45FAF16F} & \texttt{0EC4C63F} & \texttt{2C796030} & \texttt{54B14EE5} \\
\texttt{72AA6713} & \texttt{4BBF8A22} & \texttt{67475760} & \texttt{3377A4F7} & \texttt{156C8D01} & \texttt{6DA4A3D4} & \texttt{0A6249C6} & \texttt{1F8F79B5} \\

\midrule
\texttt{702D2FA9} & \texttt{3DB362C8} & \texttt{0524486F} & \texttt{33F61985} & \texttt{6B856497} & \texttt{6CA7D930} & \texttt{41DDC3E4} & \texttt{34D4A422} \\
\texttt{261B29F6} & \texttt{0206F5C8} & \texttt{0B613322} & \texttt{0C438E85} & \texttt{1E8C0351} & \texttt{794AE943} & \texttt{53124E30} & \texttt{62E2A27D} \\
\texttt{4F98B8A9} & \texttt{770F920E} & \texttt{2F7CEF1C} & \texttt{5430F397} & \texttt{17EBC5BB} & \texttt{3A91DF6F} & \texttt{46FF7E43} & \texttt{65C01FDA} \\
\texttt{7E6854E4} & \texttt{48BA050E} & \texttt{21399451} & \texttt{285E52BB} & \texttt{19AEBEF6} & \texttt{10C9781C} & \texttt{5D57357D} & \texttt{5A7588DA} \\

\midrule
\texttt{50C3CE5B} & \texttt{01D775A3} & \texttt{1A7F3E9D} & \texttt{37052449} & \texttt{08B0B349} & \texttt{3E62E2A3} & \texttt{39405F04} & \texttt{7DB9D48F} \\
\texttt{66119FB1} & \texttt{4C4938C2} & \texttt{2B8FD2D0} & \texttt{25CAA99D} & \texttt{6F76595B} & \texttt{74DE1265} & \texttt{06F5C804} & \texttt{61332216} \\
\texttt{302799EE} & \texttt{6854E4FC} & \texttt{4B6B8565} & \texttt{1318F877} & \texttt{73FCAFC2} & \texttt{7A9B6928} & \texttt{0F920EEE} & \texttt{1D5D833A} \\
\texttt{452EFE28} & \texttt{59A408B1} & \texttt{22E8143A} & \texttt{5E86B516} & \texttt{143A45D0} & \texttt{57E173FC} & \texttt{420C438F} & \texttt{2CAD6F77} \\

\midrule
\texttt{2DAE1593} & \texttt{24C9D379} & \texttt{1DDC3E48} & \texttt{669022C3} & \texttt{61B29F64} & \texttt{442D84CC} & \texttt{2A8CA834} & \texttt{58A77255} \\
\texttt{430F396B} & \texttt{745FAF17} & \texttt{1AFE83EF} & \texttt{36065EAD} & \texttt{08310E3B} & \texttt{737D12B0} & \texttt{0156C8D1} & \texttt{0F13B39C} \\
\texttt{3124E30A} & \texttt{06747576} & \texttt{384325E0} & \texttt{56E20918} & \texttt{23EB6EDE} & \texttt{68D5598E} & \texttt{7A1AD45A} & \texttt{13994505} \\
\texttt{7D3869FD} & \texttt{14BBF8A2} & \texttt{51C0B4BF} & \texttt{4D4A4226} & \texttt{6FF7E429} & \texttt{4A68FF81} & \texttt{3F619847} & \texttt{5F85CFF2} \\
\bottomrule
\end{tabular}
\end{table}

\begin{table}\centering\caption{$H_{17}$ through $H_{24}$}
\begin{tabular}{@{}cccccccc@{}}\toprule\label{table:H32_3}
\texttt{22977F14} & \texttt{1CA02F82} & \texttt{7AB1B033} & \texttt{3FCAFC2E} & \texttt{2CD20459} & \texttt{580C163C} & \texttt{01FDACB8} & \texttt{3E1D898D} \\
\texttt{74F4CB7E} & \texttt{1332216C} & \texttt{4B14EE4B} & \texttt{0E6FA256} & \texttt{59DB639F} & \texttt{1D775A21} & \texttt{3058F2C0} & \texttt{579E18D2} \\
\texttt{0FB8D7F5} & \texttt{4486E0A5} & \texttt{67EC3309} & \texttt{23400AB7} & \texttt{002AD91B} & \texttt{687E3DE7} & \texttt{12E554CF} & \texttt{2D0571FA} \\
\texttt{45519506} & \texttt{4AC39BE8} & \texttt{663B46AA} & \texttt{7523BEDD} & \texttt{69A94844} & \texttt{7B66C590} & \texttt{56496D71} & \texttt{318F8763} \\

\midrule
\texttt{25B5C2B3} & \texttt{60CE8EAE} & \texttt{7201037A} & \texttt{2A27CC5D} & \texttt{06DF111F} & \texttt{2462B710} & \texttt{377A4F67} & \texttt{427328A1} \\
\texttt{5F2EAB9B} & \texttt{2BF0B9FE} & \texttt{6E8BF5E3} & \texttt{094D1FF1} & \texttt{73D676D9} & \texttt{6F5C8040} & \texttt{1A55E786} & \texttt{36AD3AC4} \\
\texttt{089A6A52} & \texttt{1B829225} & \texttt{6119FB0D} & \texttt{50BCA575} & \texttt{516BD0D6} & \texttt{43A45D02} & \texttt{7C447837} & \texttt{070864BC} \\
\texttt{14109CCB} & \texttt{15C7E968} & \texttt{4DE1264F} & \texttt{7D930D94} & \texttt{38E84189} & \texttt{4C3653EC} & \texttt{5EF9DE38} & \texttt{393F342A} \\

\midrule
\texttt{13E62E2B} & \texttt{39BE8958} & \texttt{1DA35566} & \texttt{503D1807} & \texttt{0129A3FF} & \texttt{7A65BF74} & \texttt{2216C266} & \texttt{571FA5A0} \\
\texttt{4CB7EE9E} & \texttt{25347FC1} & \texttt{7D4702D3} & \texttt{37FBF215} & \texttt{45D02874} & \texttt{14C4938C} & \texttt{6F888F07} & \texttt{61CDF44A} \\
\texttt{2C53B92B} & \texttt{7302799E} & \texttt{0F6CD8B2} & \texttt{595ADEED} & \texttt{5E78634A} & \texttt{3E9C34FF} & \texttt{4B955339} & \texttt{060B1E58} \\
\texttt{2B71048C} & \texttt{66EF49ED} & \texttt{1A81E8C1} & \texttt{084E6515} & \texttt{68AA32A0} & \texttt{30D94FB2} & \texttt{7420C439} & \texttt{42F295D3} \\

\midrule
\texttt{028748BA} & \texttt{1048C56E} & \texttt{0E457B4D} & \texttt{49B97FEA} & \texttt{722BDA61} & \texttt{3CB0182C} & \texttt{0722BDA7} & \texttt{2E7F95F8} \\
\texttt{2BDA60E5} & \texttt{35D7DEC6} & \texttt{7EE9E996} & \texttt{40DEB900} & \texttt{778E2F7C} & \texttt{6983915F} & \texttt{30722BDB} & \texttt{3915ED31} \\
\texttt{7B4C1C8B} & \texttt{5ED30723} & \texttt{60E457B5} & \texttt{6541A2A8} & \texttt{6C266442} & \texttt{15ED3073} & \texttt{22BDA60F} & \texttt{27185312} \\
\texttt{4C1C8AF7} & \texttt{192F0384} & \texttt{0BE08E50} & \texttt{1C8AF699} & \texttt{457B4C1D} & \texttt{57B4C1C9} & \texttt{521134D4} & \texttt{5B76F23E} \\

\midrule
\texttt{4844D352} & \texttt{11B569D6} & \texttt{71513E63} & \texttt{0D3F9F4F} & \texttt{4F666EF5} & \texttt{64BC0E10} & \texttt{4601A81F} & \texttt{3D4DB494} \\
\texttt{18D2AF3C} & \texttt{21C7420D} & \texttt{2F823940} & \texttt{5DA9E321} & \texttt{639EB3B7} & \texttt{1697D471} & \texttt{7836F889} & \texttt{037AE402} \\
\texttt{5A8B5E86} & \texttt{045859A5} & \texttt{6DDBC8FA} & \texttt{6AF9755D} & \texttt{767383C4} & \texttt{412315B8} & \texttt{54CE25CB} & \texttt{1FF0129B} \\
\texttt{7F14452E} & \texttt{26E5FFAA} & \texttt{3308CFD9} & \texttt{3A6F0933} & \texttt{342A727E} & \texttt{53EC986C} & \texttt{28A084E7} & \texttt{0A1D22E8} \\

\midrule
\texttt{2F563607} & \texttt{4CE25CAB} & \texttt{39EB3B6D} & \texttt{308CFD87} & \texttt{4BC0E10C} & \texttt{3EC986CA} & \texttt{75F7B19A} & \texttt{11616691} \\
\texttt{64680157} & \texttt{0EBBAD11} & \texttt{42A727E6} & \texttt{1806A07B} & \texttt{099910B6} & \texttt{541A2A8C} & \texttt{07DC6BFB} & \texttt{7C907770} \\
\texttt{28748BA0} & \texttt{2631F0ED} & \texttt{5D7DEC66} & \texttt{1F241DDC} & \texttt{1643DB36} & \texttt{21134D4A} & \texttt{72D50C3D} & \texttt{6D0FC7BD} \\
\texttt{6A2D7A1A} & \texttt{00FED65C} & \texttt{7BB2CAD7} & \texttt{634ABCF0} & \texttt{5A5F51C1} & \texttt{45859A41} & \texttt{37AE4020} & \texttt{5338972B} \\

\midrule
\texttt{78E2F7CE} & \texttt{1B569D62} & \texttt{4FB261B2} & \texttt{41F71AFF} & \texttt{22437053} & \texttt{76A78C83} & \texttt{2561CDF4} & \texttt{5E2DD17F} \\
\texttt{67383C4E} & \texttt{1513E62F} & \texttt{2B24B6B9} & \texttt{697D4703} & \texttt{590F6CD8} & \texttt{03AEEB45} & \texttt{4890DC15} & \texttt{574A1795} \\
\texttt{46D5A758} & \texttt{34FE7D39} & \texttt{0AC92DAF} & \texttt{12315B88} & \texttt{71853124} & \texttt{2C060B1E} & \texttt{6E5FFAA4} & \texttt{0DEB9008} \\
\texttt{601A81E9} & \texttt{048C56E2} & \texttt{3D99BBD3} & \texttt{5068AA32} & \texttt{7FC04A69} & \texttt{1C7420C5} & \texttt{3ABB0674} & \texttt{33DCC09E} \\

\midrule
\texttt{075DD689} & \texttt{598ED1AA} & \texttt{332216C2} & \texttt{7F3E9C35} & \texttt{0A37FBF3} & \texttt{0918ADC4} & \texttt{7254B14F} & \texttt{5AA1879D} \\
\texttt{609B3C9B} & \texttt{3D676D8F} & \texttt{4B415C7E} & \texttt{717BE778} & \texttt{15925B5D} & \texttt{21ED9B16} & \texttt{462B7104} & \texttt{7C11CA02} \\
\texttt{6EDE47D6} & \texttt{3E483BB8} & \texttt{22C2CD21} & \texttt{54E4FCD0} & \texttt{2C87B66C} & \texttt{486E0A49} & \texttt{18F87627} & \texttt{16BD0D6A} \\
\texttt{1BD72010} & \texttt{300D40F5} & \texttt{63B46AAC} & \texttt{57CBAAE7} & \texttt{6DF111E1} & \texttt{45042733} & \texttt{047280BE} & \texttt{2FA8E05B} \\
\bottomrule
\end{tabular}
\end{table}

\begin{table}\centering\caption{$H_{25}$ through $H_{32}$}
\begin{tabular}{@{}cccccccc@{}}\toprule\label{table:H32_4}
\texttt{16E8BF5F} & \texttt{597007F6} & \texttt{483BB87C} & \texttt{712E554D} & \texttt{2FFD526E} & \texttt{5393F342} & \texttt{3A10621D} & \texttt{42D84CC8} \\
\texttt{251EA6DA} & \texttt{0D40F461} & \texttt{4C9D3785} & \texttt{3EB6EDE4} & \texttt{467EC331} & \texttt{7BCDA1F9} & \texttt{5DD6880F} & \texttt{6E20918A} \\
\texttt{6A861E73} & \texttt{57357CBB} & \texttt{124E30A6} & \texttt{64C3653E} & \texttt{07A300D5} & \texttt{7F6B2E00} & \texttt{7588DAB4} & \texttt{03058F2C} \\
\texttt{30F396A9} & \texttt{21B82923} & \texttt{1C0B4BEB} & \texttt{34551950} & \texttt{18ADC412} & \texttt{09E67B98} & \texttt{6065EAC7} & \texttt{2B5BDD97} \\

\midrule
\texttt{47A9B692} & \texttt{27CC5C55} & \texttt{20EEE1F2} & \texttt{62B71048} & \texttt{3BC717BE} & \texttt{70789D9C} & \texttt{173FCAFC} & \texttt{02D2FA8F} \\
\texttt{775A203B} & \texttt{799EE604} & \texttt{32216C26} & \texttt{4E4FCD0A} & \texttt{29089A6A} & \texttt{6C73D677} & \texttt{524486E1} & \texttt{55663B46} \\
\texttt{0C163CB0} & \texttt{05F04728} & \texttt{101D775B} & \texttt{19FB0CC3} & \texttt{408B0B35} & \texttt{0B348117} & \texttt{6595ADEF} & \texttt{496D70AD} \\
\texttt{3503D181} & \texttt{5BA2FD79} & \texttt{7EBC5BA3} & \texttt{3CE5AA19} & \texttt{6B516BD0} & \texttt{2E2A27CD} & \texttt{1ED9B164} & \texttt{5C8040DE} \\

\midrule
\texttt{58737D12} & \texttt{680156C9} & \texttt{4FE7D387} & \texttt{7F95F85C} & \texttt{38972AA7} & \texttt{7DEC66BA} & \texttt{1F71AFE9} & \texttt{6A78C82F} \\
\texttt{4D9E4D61} & \texttt{2F038432} & \texttt{544F98B9} & \texttt{71D08311} & \texttt{3AEEB441} & \texttt{41A2A8CA} & \texttt{06A07A31} & \texttt{0A9C9F9A} \\
\texttt{134D4A42} & \texttt{1134D4A4} & \texttt{34ABCF0C} & \texttt{04D9E4D7} & \texttt{5636065F} & \texttt{2146FF7F} & \texttt{1D08310F} & \texttt{66442D84} \\
\texttt{73A91DF7} & \texttt{36D251EA} & \texttt{43DB362C} & \texttt{08E5017C} & \texttt{233F6199} & \texttt{5A0AE3F4} & \texttt{643DB362} & \texttt{2D7A1AD4} \\

\midrule
\texttt{2EAB9ABF} & \texttt{109CCA29} & \texttt{274DE127} & \texttt{791F5B76} & \texttt{4ECE7078} & \texttt{5C01FDAC} & \texttt{70F920EE} & \texttt{32A0D154} \\
\texttt{55E78634} & \texttt{6CF26B05} & \texttt{0571FA5A} & \texttt{77DB9D49} & \texttt{49ECCDDF} & \texttt{47280BE0} & \texttt{1E580C16} & \texttt{6BD0D6A2} \\
\texttt{206F5C80} & \texttt{3B46AACC} & \texttt{5B23400B} & \texttt{6236AD3A} & \texttt{7E3DE6D1} & \texttt{0C9781C2} & \texttt{6514109D} & \texttt{025347FD} \\
\texttt{29892718} & \texttt{52C53B93} & \texttt{197AB1B1} & \texttt{17BE778E} & \texttt{0BB53C65} & \texttt{400AB647} & \texttt{3C64176B} & \texttt{35826CF3} \\

\midrule
\texttt{4947A9B6} & \texttt{29A3FE03} & \texttt{10B61332} & \texttt{6BFA0FB9} & \texttt{27E6854E} & \texttt{5BDD9657} & \texttt{320BB53D} & \texttt{1EF3687F} \\
\texttt{20918ADC} & \texttt{022C2CD3} & \texttt{4075DD69} & \texttt{62C87B66} & \texttt{5CAA99C5} & \texttt{3C4ECE70} & \texttt{65BF74F4} & \texttt{5598ED1A} \\
\texttt{70524487} & \texttt{79603058} & \texttt{7E173FCA} & \texttt{52EFE288} & \texttt{198467ED} & \texttt{3B39C1E2} & \texttt{0B1E580C} & \texttt{6C8D002B} \\
\texttt{357CBAAF} & \texttt{77254B15} & \texttt{055B2341} & \texttt{2ED4F191} & \texttt{17C11CA0} & \texttt{4702D2FB} & \texttt{0C69579E} & \texttt{4E30A624} \\

\midrule
\texttt{0EEE1F24} & \texttt{239405F0} & \texttt{1C2192F0} & \texttt{6E0A4891} & \texttt{1546541A} & \texttt{4A1794AF} & \texttt{7BE778E2} & \texttt{75A203AF} \\
\texttt{51BFDF91} & \texttt{0789D9CE} & \texttt{3F1EF369} & \texttt{2AF3C31A} & \texttt{315B8824} & \texttt{4452EFE2} & \texttt{6928F536} & \texttt{676D8E7B} \\
\texttt{7CC5C545} & \texttt{09CCA283} & \texttt{2DD17EBD} & \texttt{7280BE08} & \texttt{58D8197B} & \texttt{43705245} & \texttt{4D352908} & \texttt{24B6B857} \\
\texttt{604F33DC} & \texttt{1B032F57} & \texttt{1264E9BD} & \texttt{36793583} & \texttt{00AB6469} & \texttt{5FFAA4DC} & \texttt{383C4ECE} & \texttt{569D6236} \\

\midrule
\texttt{5A203AEF} & \texttt{163CB018} & \texttt{64E9BC25} & \texttt{1879CB55} & \texttt{33A3ABB0} & \texttt{6AACC768} & \texttt{71048C56} & \texttt{0D6A2D7A} \\
\texttt{5ADEECB3} & \texttt{7FBF2147} & \texttt{0D94FB26} & \texttt{7F41F71B} & \texttt{4FCD0A9C} & \texttt{280BE08E} & \texttt{71FA5A0A} & \texttt{4F33DCC0} \\
\texttt{18871D09} & \texttt{26B04D9F} & \texttt{28F536D2} & \texttt{546541A2} & \texttt{03D1806B} & \texttt{549B97FE} & \texttt{418871D1} & \texttt{032F5637} \\
\texttt{335D7DEC} & \texttt{3D1806A1} & \texttt{264E9BC3} & \texttt{16C26644} & \texttt{64176A79} & \texttt{4176A78D} & \texttt{6A521134} & \texttt{3DE6D0FD} \\

\midrule
\texttt{24370525} & \texttt{4A9629DD} & \texttt{08CFD867} & \texttt{3EE35FD1} & \texttt{50427329} & \texttt{561CDF44} & \texttt{3784993B} & \texttt{01A81E8D} \\
\texttt{2D50C3CF} & \texttt{69579E18} & \texttt{75DD6881} & \texttt{7CBAAE6B} & \texttt{121B8293} & \texttt{7383C4EC} & \texttt{0E91740A} & \texttt{1B7C4479} \\
\texttt{07F6B2E0} & \texttt{31DA3556} & \texttt{1D22E814} & \texttt{4CC885B0} & \texttt{6F093275} & \texttt{666EF49F} & \texttt{5925B5C3} & \texttt{38BDF3BC} \\
\texttt{45AF435A} & \texttt{14452EFE} & \texttt{603058F2} & \texttt{43F1EF37} & \texttt{2269A948} & \texttt{5F7B19AE} & \texttt{7AE40206} & \texttt{2B0E6FA2} \\
\bottomrule
\end{tabular}
\end{table}

\bibliographystyle{plain}
\bibliography{uleth-refs}

\end{document}